\documentclass[journal,twoside,web]{ieeecolor}
\usepackage{generic}
\usepackage{cite}
\usepackage{amsmath,amssymb,amsfonts}
\usepackage{algorithmic}
\usepackage{graphicx}
\usepackage{algorithm,algorithmic}{\tiny }
\usepackage{hyperref}
\hypersetup{hidelinks=true}
\usepackage{textcomp}
\newtheorem{lemma}{Lemma}
\usepackage{verbatim}
\usepackage{url}
\usepackage[caption=false,font=normalsize,labelfont=sf,textfont=sf]{subfig}
\usepackage{algorithm}
\usepackage{array}
\usepackage{color}
\def\BibTeX{{\rm B\kern-.05em{\sc i\kern-.025em b}\kern-.08em
		T\kern-.1667em\lower.7ex\hbox{E}\kern-.125emX}}
\begin{document}

\title{Observer-Based Event-Triggered Secure Consensus Control for Multi-Agent Systems}

\author{Jingyao~Wang, Zeqin~Zeng, Jinghua~Guo, and Zhisheng~Duan 
\thanks{This work was supported by the National Nature Science Foundation of China under Grant Nos. $62473323$, $52372419$, $62173006$, T$2121002$, and $61803319$. }
\thanks{J. Wang and Z. Zeng are with the School of Aerospace Engineering, Xiamen University, Xiamen 361005, P. R. China (e-mail: wangjingyao1@xmu.edu.cn; zzq19991108@163.com).}
\thanks{J. Guo is with the Pen-Tung Sah Institute of Micro-Nano Science and Technology, and the Institute of Artiffcial Intelligence, Xiamen University, Xiamen 361005, P. R. China (e-mail: guojing$\_$0701@live.cn). \textit{(Corresponding author: Jinghua Guo.)}}
\thanks{Z. Duan is with the Department of Automation, Xiamen University, Xiamen 361000, P. R. China, and also with the Department of Mechanics and Engineering Science, Peking University, Beijing 100871, China (e-mail: duanzs@pku.edu.cn).}}%}

%\markboth{}%
%{}


\maketitle

\begin{abstract}
This study delves into the intricate challenges encountered by multi-agent systems (MASs) operating within environments that are subject to deception attacks and Markovian randomly switching topologies, particularly in the context of event-triggered secure consensus control. To address these complexities, a novel observer-based distributed event-triggered control scheme is introduced. This approach uses local information to dynamically adjust its triggered conditions, thereby enhancing the utilization of network resources. Additionally, the design of the observer based secure consensus controller is distributed, leveraging the local information of each individual agent. Furthermore, our event-triggered mechanism theoretically precludes the occurrence of Zeno behavior in the triggering time series. Finally, simulation results underscore the superiority of our proposed method when compared to existing techniques, thereby validating its effectiveness and applicability in the event-triggered secure consensus control of MASs.
\end{abstract}

\begin{IEEEkeywords}
 Multi-agent systems, deception attacks, randomly switching topologies, event-triggered, secure control
\end{IEEEkeywords}

\section{Introduction}\label{sec:introduction}
\IEEEPARstart{T}{he} network security holds paramount importance in the distributed control of multi-agent systems (MASs), comprising agents that exchange information through an interactive network. However, achieving a fully secure network environment remains elusive, as attackers can easily capitalize on the data transmitted among agents. Unlike naturally occurring disturbances \cite{Sun2022} and faults \cite{Chen2015} in traditional systems, cyberattacks are deliberately and unpredictably perpetrated by adversaries, aiming to compromise the overall performance of the system and potentially induce instability. Given the intricate nature and inherent vulnerabilities of distributed MASs, it is crucial to devise robust security measures that effectively mitigate the potential risks posed by such attacks. Therefore, a thorough investigation into the necessary security mechanisms for maintaining consensus among agents in MASs is imperative.

Deception attacks consist of unauthorized access perpetrated by attackers, subsequently leading to the manipulation of measured data and control signals. The inherent stealthiness of these attacks poses significant challenges in their detection, potentially causing substantial degradation in system performance. The attack mechanism investigated in this study entails the injection of falsified data into the information channel, ultimately culminating in misleading estimations within the MASs. The fault-tolerant control issue for MASs facing deception attacks is tackled in \cite{WOS17}, introducing a distributed impulsive controller as a solution. In \cite{WOS18}, the bipartite synchronization challenge under deception attacks is addressed by developing an impulsive control scheme. Zhang \textit{et al.} introduced in \cite{WOS19} an observer-based non-fragile $H_{\infty}$ consensus controller tailored for individual agents. This controller ensures the preservation of the designated disturbance attenuation level, even in the face of deception attacks. 

The preceding attack analysis assumes a static topology, overlooking the dynamic nature of agent-to-agent connections in real-world engineering scenarios. In these applications, various factors—including signal power attenuation, thermal noise \cite{WOS7}, environmental obstacles, and constraints on sensor communication ranges—result in temporal variations in the connections between agents. To alleviate the stringent requirements on the connection graph, Markovian randomly switching communication graphs have garnered significant attention in the research community. Li \textit{{et al.}} \cite{WOS20} and Wang \textit{{et al.}} \cite{WOS21} focused on addressing the mean square consensus issue for continuous-time MASs, considering the presence of Markovian randomly switching topologies. To tackle the complexities arising from the interplay between heterogeneous dynamics and input saturation constraints, Wang \textit{{et al.}}\cite{WOS23} introduced distributed adaptive observers and antiwindup controllers, aimed at achieving output consensus for heterogeneous MASs under Markovian randomly switching topologies.

In the face of network attacks, traditional continuous communication and control strategies often appear inadequate. They not only exacerbate the communication burden but also potentially increase the risk of being attacked due to continuous exposure to the network environment. Additionally, the limited energy supply of agents and the constraints of network bandwidth are practical factors that must be taken into account. Conversely, the event-triggered mechanism given in\cite{WOS8,WOS9,WOS10,WOS11,WOS12} can not only significantly reduce the frequency of communication and data transmission, alleviating the communication burden, but also improve the response speed and security of the system as shown in \cite{Han2023}.

It is crucial to develop an event-triggered secure controller that can effectively mitigate deception attacks and adapt to randomly switching topologies encountered by MASs. Such a controller would not only ensure the resilience and security of MASs operating in dynamic and adversarial environments, but also optimize communication resource utilization. This paper aims to develop such an event-triggered secure consensus control scheme. The contributions of this paper are as follows:
\begin{itemize}
	\item [(1)]This article proposes a novel event-triggered mechanism that leverages attacked output information to minimize unnecessary communication and thereby enhance the efficient utilization of system resources.
	\item [(2)]The proposed controller facilitates the attainment of secure consensus among agents within Markovian randomly switching communication topologies, where the communication graph may become disconnected. Consequently, it relaxes the conventional assumption that the undirected graph must be fully connected, as stipulated in \cite{WOS13,graph2,graph3}. This advancement affords a more adaptable approach to the consideration of communication topology, thus alleviating traditional constraints. 
	\item [(3)]To the best of our knowledge, this is the first time that a distributed event-triggered consensus mechanism has been designed for agents under deception attacks and randomly switching communication topologies, utilizing solely the local information available to each agent. This approach represents a novel contribution to the published literature in this area. 
\end{itemize}

The remainder of this article is organized as follows. In Section II, we introduce the basic notation from graph theory and formulate the problem. Section III presents the observer-based distributed  event-triggered secure control scheme for MASs. For further clarification, numerical examples are provided in Section IV. Finally, Section V concludes the article with remarks.

\section{PROBLEM FORMULATION}
In this section, we introduce the notation from graph theory and the models of system dynamics and deception attack that are crucial for the problem formulation.
\subsection{Notation}
Index the nodes from $1$ to $N$, with the $i$-th agent denoted as agent $i$. Let ${\cal V}$ represent the set of all nodes, defined as ${\cal V} = \left\{ 1,2,...,N \right\}$.  The set of edges, denoted by ${\cal E}$, consists of all the connections between nodes. The communication graph is denoted by ${\cal G} = \left( {{\cal V},{\cal E}} \right)$. The adjacency matrix of G is represented by ${\cal A} = \left[ {{a_{ij}}} \right] \in {{\rm{R}}^{N \times N}}$, where ${a_{ij}} = 1$ if there is an edge between $i$ and $j$ and ${a_{ij}} = 0$ otherwise. The Laplacian matrix of ${\cal G}$ is denoted by ${\cal L} = {\left( {{l_{ij}}} \right)_{N \times N}}$, where ${{l_{ij}}}$ = $ - {a_{ij}}$ for $i \ne j$ and ${l_{ii}}$ = $\sum\nolimits_{j = 1}^N {{a_{ij}}} $.

At any given time $t$, let ${\cal G}\left( t \right)$ represent the pair $\left( {{\cal V},{\cal E}\left( t \right)} \right)$, where the edge set ${{\cal E}\left( t \right)}$ evolves over time. If the changes in ${\cal G}\left( t \right)$ are governed by a random process, then the communication graphs are said to undergo randomly switching.

The continuous-time Markov process is represented by $\left\{ {\sigma \left( t \right),t \in {{\rm{R}}_ + }} \right\}$. Its infinitesimal generator matrix is denoted by  $\Upsilon = \left[ {{w_{ij}}} \right]$. The transition probability from state $i$ to state $j$ during a small time interval ${\bar \tau }$ is given by ${w_{ij}}\bar \tau  + o\left( {\bar \tau } \right)$, where $o\left( {\bar \tau } \right)$ satisfies $\mathop {\lim }\nolimits_{\bar \tau  \to 0} \frac{{o\left( {\bar \tau } \right)}}{{\bar \tau }} = 0$. If the process remains in state $i$, the transition probability is  $1 + {w_{ii}}\bar \tau  + o\left( {\bar \tau } \right)$. The row summation of the transition rate matrix $\Upsilon $, is zero, i.e., $\Upsilon \textbf{1} = 0$.

%\section{PROBLEM FORMULATION}
\subsection{System model}
Consider a MAS comprising $N$ agents, with the dynamics of the $i$-th agent described by
\begin{IEEEeqnarray}{l} \label{eq2}
	{{\dot x}_i}\left( t \right) = A{x_i}\left( t \right) + B{u_i}\left( t \right),i = 1,2,...,N, \nonumber\\
	{y_i}\left( t \right) = C{x_i}\left( t \right),
\end{IEEEeqnarray}
in which ${x_i}\left( t \right) \in {{\rm{R}}^n}$ represents the state, ${u_i}\left( t \right) \in {{\rm{R}}^m}$ denotes the control input, and ${y_i}\left( t \right) \in {{\rm{R}}^p}$ corresponds to the measured output of the $i$-th agent. The matrices $A$, $B$, and $C$ have rational dimensions.

The connection graph among $N$ agents ${\cal G}\left( t \right)$ randomly alternates among $s$ distinct graphs, denoted as ${\cal G}\left( t \right) \in \left\{ {{{\cal G}_1},{{\cal G}_2},...,{{\cal G}_s}} \right\}$. Specifically, ${\cal G}\left( t \right)$ equals ${{\cal G}_p}$ precisely when $\sigma \left( t \right)$ belongs to the set $\left\{ {1,2,...,s} \right\}$. Each of these $s$ graphs is represented as ${{\cal G}_1} = \left( {{\cal V},{{\cal E}_1}} \right),...,{{\cal G}_s} = \left( {{\cal V},{{\cal E}_s}} \right)$. The union graph, denoted as ${\cal G}_{\text{un}}$, is characterized as the graph whose vertex set is ${\cal V}$ and whose edge set is the union of all edge sets ${\cal E}_p$ for $p$ ranging from $1$ to $s$. Furthermore, the corresponding adjacency matrix and Laplacian matrix of ${{\cal G}_{\text{un}}}$ are denoted by $\hat {\cal A} = \left[ {{{\hat a}_{ij}}} \right]$ and $\hat {\cal L} = [ {{{\hat l}_{ij}}} ]$, respectively.
\newtheorem{assumption}{Assumption}[section]
\begin{assumption}\label{assumption1}
	The union graph ${{\cal G}_{\text{un}}}$ of the $s$ undirected graphs is undirected and connected.
\end{assumption}

\begin{assumption}\label{assumption2}
The communication graph ${\cal G}(t)$ stochastically alternates among $s$ undirected graphs, and this switching procedure of these connection graphs is guided by a continuous-time ergodic Markov process possessing a transition rate matrix $\Upsilon$.
\end{assumption}

\newtheorem{remark}{Remark}
From Assumption \ref{assumption2}, we can deduce that the Markov process possesses a unique invariant distribution denoted as $\Pi = {\left[ {{\Pi _1},{\Pi _2},...,{\Pi _s}} \right]^T}$. This distribution satisfies the conditions that ${\Pi ^T}\Upsilon  = 0$ and ${\Pi ^T}\textbf{1} = 1$ while maintaining non-negative values for ${\Pi _p}$  with $p$ ranging from $1$ to $s$. Furthermore, this invariant distribution remains constant at $\Pi $ for all $t \ge 0$.

\begin{remark}
The communication graphs considered in this manuscript can be disconnected, which is in accordance with the information connection graph in practice.	\end{remark}
\subsection{Deception attacks}
The output signals under deception attacks of agent $i$ are modeled as
\begin{equation}\label{eq4}
	{\mathord{\buildrel{\lower3pt\hbox{$\scriptscriptstyle\smile$}} 
			\over y} _i}\left( t \right) = {y_i}\left( t \right) + {\alpha _i}\left( t \right){\varepsilon _i}\left( t \right),
\end{equation}
where ${\varepsilon _i}\left( t \right)$ represent the attack signals injected into the transmission channels.  The random variables ${\alpha _i}\left( t \right)$ indicate the presence of attacks on agent $i$'s transmission channels, where ${\alpha _i}\left( t \right) = 1$ indicates an attack and ${\alpha _i}\left( t \right) = 0$ indicates no attack. These random variables ${\alpha _i}\left( t \right)$ follow a Bernoulli distribution, with the probability of an attack given by
\begin{equation}\label{eq5}
	\text{prob}\left\{ {{\alpha _i}\left( t \right) = 1} \right\} = {{\bar \alpha }_i}
\end{equation}
where the probability distribution ${{\bar \alpha }_i}$, $i = 1,2,...,N$, are known, rewritten as $F =$ diag$\left\{ {{{\bar \alpha }_1},{{\bar \alpha }_2},...,{{\bar \alpha }_N}} \right\}$, in which diag$\left\lbrace  {\cdot} \right\rbrace  $ stands for the diagonal matrix.

\begin{assumption}\label{assumption3}
	Deception attack signals ${\varepsilon _i}(t)$ are presumed to fulfill the subsequent criterion:
	\begin{equation}\label{eq6}
		{\left\| {{\varepsilon}\left( t \right)} \right\|^2} \le \tau, 
	\end{equation}
	where $\varepsilon \left( t \right) = \left[  {{\varepsilon _1}^T\left( t \right),{\varepsilon _2}^T\left( t \right),...,{\varepsilon _N}^T\left( t \right)} \right]^T$ and $\tau $ is a known positive constant.
\end{assumption}

 The objective of the manuscript is to design a controller, based on the event-triggered mechanisms, for the MAS with dynamics given in (\ref{eq2}) under deception attacks and Markovian randomly switching topologies such that the sate of each agent can converges to a small neighborhood of the average state of all agents in the mean square sense.

\begin{figure}[!t]
	\centerline{\includegraphics[width=\columnwidth]{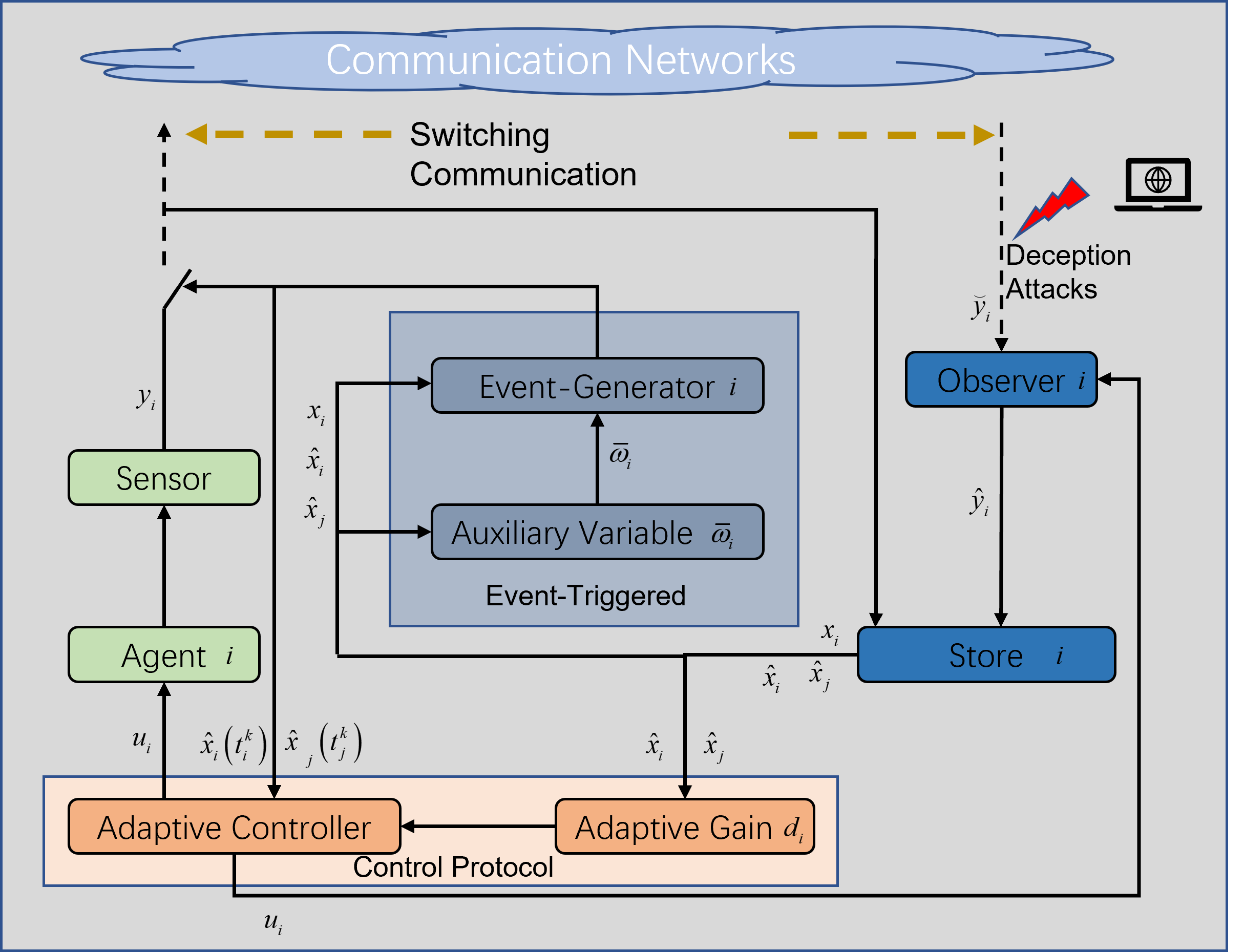}}
	\caption{Framework for the distributed secure control protocol of the MASs under deception attacks and Markovian randomly switching communication topologies.}
	\label{frame}
\end{figure}

\section{MAIN RESULTS}
In this section, we proceed to design the distributed secure consensus control protocol for the MAS under deception attacks and randomly switching communication topologies. The framework for the secure consensus control protocol is shown in Fig.\ref{frame}, where event-triggered mechanism is introduced for updating control signals at specific event-triggered instants.  

Considering that in practice only the output information is available, the following observer is introduced to estimate the state information
\begin{IEEEeqnarray}{L} \label{eq3}
	{{\dot{\hat x}}_i}\left( t \right) = A{{\hat x}_i}\left( t \right) + B{u_i}\left( t \right) + G\left( {{{\mathord{\buildrel{\lower3pt\hbox{$\scriptscriptstyle\smile$}} 
					\over y} }_i}\left( t \right) - {{\hat y}_i}\left( t \right)} \right), \nonumber\\
	{{\hat y}_i}\left( t \right) = C{{\hat x}_i}\left( t \right),i=1,\ldots,N,
\end{IEEEeqnarray}
where ${{\hat x}_i}\left( t \right) \in {{\rm{R}}^n}$ and ${{\hat y}_i}\left( t \right) \in {{\rm{R}}^p}$ are the state and output of the observer system for the $i$-th agent, respectively. $G \in {{\rm{R}}^{n \times p}}$ is the observer gain to be determined.

Given agent $i=1,\ldots,N$, its network induced state estimation error is defined as
\begin{equation}
	{m_i}(t) = {\hat x_i}(t_k^i) - {\hat x_i}(t),
\end{equation}
where ${t_k^i}$ is the latest event-triggered instant and ${\hat x_i}(t_k^i)$ denotes the estimated state of agent $i$ at time ${t_k^i}$. The determination of the triggered time sequence $\{{t_k^i}, k\in \mathbb{N}^+\}$ is governed by an event generator that utilizes the following conditions
\begin{align} \label{eq9}
	t_1^i &= 0, \nonumber\\
	t_{k + 1}^i &= \inf\limits_{l > t_k^i} \bigg\{ {l:{\iota _i}( {m_i^T(t)\Gamma {m_i}(t) - {o_i}{\upsilon _i}\tilde \xi _{\rm{i}}^T(t)\Gamma {{\tilde \xi }_i}(t)} )} \nonumber\\
	&   \ge {\varpi _i}(t),\forall t \in ( {t_k^i,l} ] \bigg\},
\end{align}
where the relative state error between neighboring agents at the event-triggered instant is denoted by ${{\tilde \xi }_i}\left( t \right) = {\xi _i}\left( {t_k^i} \right)$, for $t=t_k^i$; ${\xi _i}(t_k^i) = \mathop \sum \limits_{j = 1}^N {a_{ij}}\left( {\sigma \left( {t_k^i} \right)} \right)\left( {{{\hat x}_j}(t_k^j) - {{\hat x}_i}(t_k^i)} \right)$; $\Gamma  = PB{B^T}P$, in which $P$ is a positive-definite symmetric matrix to be determined; ${o_i} \in \left( {0,1} \right)$; ${\upsilon _i}$ and ${{\iota _i}}$ are positive constants; the threshold parameter ${\varpi _i}(t)$ is an auxiliary variable satisfying
\begin{align}\label{auxiliary variable}
	{\dot \varpi _i}(t) =&  - {\eta _i}{\varpi _i}(t) \nonumber\\
	&+ {\varsigma _i}\left[ {{o_i}{\upsilon _i}\tilde \xi _{\rm{i}}^T\left( t \right)\Gamma {{\tilde \xi }_i}\left( t \right) - m_i^T(t)\Gamma {m_i}(t)} \right],
\end{align}
with ${\varpi _i}(0) > 0$, ${\eta _i} > 0$, ${\varsigma _i} > 0$.

\begin{lemma}
The threshold parameters designed by (\ref{auxiliary variable}) are positive, i.e., ${{\varpi _i}\left( t \right)}>0, i=1,\ldots,N$ for $t\geq 0$.
\end{lemma}
\textit{Proof}:
Consider two consecutive event-triggered instants $t_k^i$ and $t_{k+1}^i$, which implies that no event is triggered for $t\in (t_k^i,\ t_{k+1}^i)$. Thus, we can get that the inequality in (\ref{eq9}) is not satisfied for $t\in (t_k^i,\ t_{k+1}^i)$. That is to say
\begin{equation}
{\iota _i}\left( {m_i^T\left( t \right)\Gamma {m_i}\left( t \right) - {o_i}{\upsilon _i}\tilde \xi _{\rm{i}}^T\left( t \right)\Gamma {{\tilde \xi }_i}\left( t \right)} \right)< {\varpi _i}\left( t \right).
\end{equation}
Substituting the above inequality into (\ref{auxiliary variable}) yields
\begin{equation}
{\dot \varpi _i}(t)>- {\eta _i}{\varpi _i}(t)-\frac{\varsigma _i}{\iota _i}{\varpi _i}(t).
\end{equation}

With some mathematical calculations, we have
 \begin{align}
{ \varpi _i}(t) \ge { \varpi _i}\left( {t_k^i} \right){e^{ - \left( {{\eta _i} + \frac{{{\varsigma _i}}}{{{\iota _i}}}} \right)\left( {t - t_k^i} \right)}},t \in \left[ {t_k^i,t_{k + 1}^i} \right). \nonumber
 \end{align}
And noting that $t_1^i = 0$, we have
\begin{align}
	{{ \varpi }_i}\left( t \right) \ge& {\varpi _i}\left( {t_k^i} \right){e^{ - \left( {{\eta _i} + \frac{{{\varsigma _i}}}{{{\iota _i}}}} \right)\left( {t - t_k^i} \right)}} \nonumber\\
	\ge& {\varpi _i}\left( {t_{k - 1}^i} \right){e^{ - \left( {{\eta _i} + \frac{{{\varsigma _i}}}{{{\iota _i}}}} \right)\left( {t - t_{k - 1}^i} \right)}}\nonumber\\
	\ge&  \cdots  \ge {\varpi _i}\left( 0 \right){e^{ - \left( {{\eta _i} + \frac{{{\varsigma _i}}}{{{\iota _i}}}} \right)t}}>0.\nonumber
\end{align}
Thus, we get that ${{\varpi _i}\left( t \right)}>0$ for $t\geq 0$.
\begin{remark}
The dynamic event-triggered strategy (\ref{eq9}) determines the triggered instant solely leveraging the local information available to each agent. This distributed communication mechanism minimizes unnecessary information communication, thereby optimizing the utilization of network communication resources.
\end{remark}

Within the proposed framework, an observer-based distributed event-triggered secure controller is introduced to solve the secure consensus control problem. The controller's formulation is given by
\begin{eqnarray}\label{controller}
	{u_i}(t) &=& {d_i}(t)K\mathop \sum \limits_{j = 1}^N {a_{ij}}\left( {\sigma \left( t \right)} \right)\left( {{{\hat x}_j}(t_k^j) - {{\hat x}_i}(t_k^i)} \right),\nonumber\\
	&& t\in \left[ {t_k^i,t_{k + 1}^i} \right),\ i=1,\ldots,N,
\end{eqnarray}
where $K$ represents the gain matrix to be determined, while ${d_i}(t)$ signifies the adaptive coupling strength that follows a specific update law  outlined as follows
\begin{eqnarray}\label{eq11}
	{\dot d_i}(t) = \left\{ \begin{array}{ll}
		{\beta _i}{\tilde \xi _{\rm{i}}^T\left( t \right)\Gamma {{\tilde \xi }_i}\left( t \right)} & \textrm{if ${d_i}(t) < {{\bar d}_i}$},\\
		0 & \textrm{if ${d_i}(t) \ge {{\bar d}_i}$},
	\end{array} \right.
\end{eqnarray}
where ${\beta _i}$ and ${\bar d_i}$ are predefined positive constants. And, the initial condition for ${d_i}(t)$ lies within the range $1 < {d_i}({0}) < {\bar d_i}$. This ensures that the coupling strength ${d_i}(t)$ satisfies ${d_i}({0}) < {d_i}({t}) < {\bar d_i}$.
In addition, to simplify the notation, the variable $t$ will be omitted from the subsequent discussion.

 Based on equations (\ref{eq2}), (\ref{eq3}), (\ref{controller}), and the definitions of ${\varepsilon _i}$ and ${{\tilde \xi }_i}$, the MAS (\ref{eq2}) and observer (\ref{eq3}) can be rewritten respectively as
\begin{align} \label{eq12}
	\dot x &= \left( {{I_N} \otimes A} \right)x + \left( {d \otimes BK} \right){\tilde \xi }, \nonumber\\
	\dot{\hat x} &=\left( {{I_N} \otimes A} \right)\hat x + \left( {d \otimes BK} \right)\tilde \xi \nonumber\\
	&\quad \quad - \left( {{I_N} \otimes GC} \right)\left( {\hat x - x} \right) + \left( {\alpha  \otimes G} \right)\varepsilon,
\end{align}
where $x = \left[  {x_1^T,...,x_N^T} \right]^T$, $\hat x = \left[  {\hat x_1^T,...,\hat x_N^T} \right]^T$, $\tilde \xi = \left[  {\tilde \xi_1^T,...,\tilde \xi_N^T} \right]^T$, $d =$ diag$\left\{ {d_1^T,...,d_N^T} \right\}$ and $\alpha=$ diag$\left\{ {\alpha_1,...,\alpha_N} \right\}$.

The state estimation error of agent $i$ is defined as ${q_i} = {\hat x_i} - {x_i}, i=1,\ldots,N$. Let $q =\left[  {q_1^T,...,q_N^T} \right]^T$. Then, based on (\ref{eq12}), the dynamics of the state estimation error is given by
\begin{equation}\label{eq14}
	\dot q = \left[ {{I_N} \otimes \left( {A - GC} \right)} \right]q + \left( {\alpha \otimes G} \right)\varepsilon.
\end{equation}

To define the average state of all agents, we introduce the notation $\bar x(t) = 1/N\sum_{i = 1}^N {{x_i}(t)} $. Subsequently, we introduce the error signal ${\delta _i} (t)= {x_i}(t) - \bar x(t)$ to represent the deviation of agent $i$'s state from the average state, where $i$ ranges from $1$ to $N$. By stacking the error signals in a column vector, we obtain $\delta(t) = \left[  {\delta _1^T(t),\cdot  \cdot  \cdot ,\delta _N^T(t)} \right]^T$. It can be shown that $\delta(t)$ can be stated as $\delta(t) = \left( {M \otimes {I_n}} \right)x(t)$, where $M = {I_N} - \frac{1}{N}{\textbf{11}^T}$.
The consensus error system can be described by the following equation
\begin{equation}\label{eq15}
	\dot \delta = \left( {{I_N} \otimes A} \right)\delta + \left( {Md \otimes BK} \right){\tilde \xi }.
\end{equation}

Now, based on the above definition, the secure consensus control problem of the MAS with dynamics described by (\ref{eq2}) can be converted to the stability problem of the associated closed-loop error systems (\ref{eq14}) and (\ref{eq15}). The following theorem presents a sufficient condition which ensures the ultimate boundedness of the consensus error.
\newtheorem{theorem}{\bf Theorem}
\begin{theorem}\label{thm1}
	Under Assumptions \ref{assumption1}-\ref{assumption3}, the  consensus error system (\ref{eq15}) exponentially converges to the set ${\bar \ell }$ in the mean-square sense with
	\begin{equation}
		\bar \ell  = \left\{ {\delta  \in {{\rm{R}}^N}|\mathbb{E}\left\{ {\left\| {\delta } \right\|} \right\} \le \sqrt {\frac{\tau }{{\chi \kappa }}} } \right\},
	\end{equation}
	if there exist positive scalars $c$, ${\eta _i}$, ${\iota _i}$, ${\varsigma _i}$, ${\rho _i}$, ${o_i} \in \left( {0,1} \right)$,  $\kappa  $, $\chi  $ and positive definite matrices $P$ and $Q$ so that the following formulas are satisfied
	\begin{IEEEeqnarray}{c} \label{eq17}
		\Lambda  + \kappa \tilde \Theta  < 0\\ \label{eq18}
		\chi I - \hat \Theta  < 0, \\
		{\rho _i} - {\varsigma _i} > 0, \\
		{\eta _i} - \frac{{{\rho _i} - {\varsigma _i}}}{{{\iota _i}}} > 0, \\
		{{\bar d}_i} > 4c + {o_i}+1,\\
		{\upsilon _i} \geq \frac{1}{{{\rho}}},
	\end{IEEEeqnarray}
	where 
	\begin{IEEEeqnarray}{c}
		\Lambda  = \left( {\begin{array}{*{20}{c}}
				{\bar \Pi \left( {\left( {{A^T}P + PA} \right) - \frac{c}{s}{\lambda _2}( {\hat {\cal L}} )\Gamma } \right)}&0\\
				0&{ \Psi }
		\end{array}} \right),\nonumber\\
		{ \Psi }  =\frac{1}{s}  \left( {{A^T}Q + QA - QGC - {{\left( {QGC} \right)}^T}} \right.\nonumber\\
		\left. { + {\lambda _M}( {F{F^T}} )QG{G^T}{Q^T} + s\left( {2c + {{\bar d}^2}} \right)\mathord{\buildrel{\lower3pt\hbox{$\scriptscriptstyle\smile$}} 
				\over \Pi }\lambda _M^2( {\hat {\cal L}} )\Gamma } \right) ,\nonumber\\
		\tilde \Theta {\rm{ = }}\left( {\begin{array}{*{20}{c}}
				P&0\\
				0&Q
		\end{array}} \right) \nonumber ,
\hat \Theta  = \left( {\begin{array}{*{20}{c}}
		{{\lambda _2}( {\hat {\cal L}} )P}&0\\
		0&Q
\end{array}} \right),
	\end{IEEEeqnarray}
	with ${{\lambda _2}( {\hat {\cal L}} )}$ be the second smallest eigenvalue of $\hat {\cal L} = \sum\nolimits_{p = 1}^s {{\cal L}\left( p \right)} $, $\mathord{\buildrel{\lower3pt\hbox{$\scriptscriptstyle\smile$}} 
		\over \Pi }  = {\max _{p = 1,2,...,s}}{\Pi _p}$, $\bar \Pi  = {\min _{p = 1,2,...,s}}{\Pi _p}$, $\rho  = \tilde c\mathord{\buildrel{\lower3pt\hbox{$\scriptscriptstyle\smile$}} 
		\over \Pi } N\left( {{N^2} + N} \right)$ and $\tilde c = \mathop {\max }_{i = 1,2,...,N} \left\{ {{{\bar d}_i} + 4c} \right\}$; and the controller gain matrix is calculated with $K = {B^T}P$.
\end{theorem}

%\begin{IEEEproof}

\textit{Proof}: Construct a candidate Lyapunov function
\begin{equation}
	V = {V_1} + {V_2} = \sum\limits_{p = 1}^s {{V_{1p}}}  + \sum\limits_{p = 1}^s {{V_{2p}}} 
\end{equation}
where ${V_{1p}} = \mathbb{E}\left[( {{Z^T}{\Theta_p} Z){\textbf{1}_{\left\{ {\sigma \left( t \right) = p} \right\}}}} \right]$, $Z= {\left[ {{\delta ^T},{q^T}} \right]^T}$, ${\Theta_p}  = \left( {\begin{array}{*{20}{c}}
		{{\cal L}\left( p \right) \otimes P}&0\\
		0&{\frac{1}{s}{I_N} \otimes Q}
\end{array}} \right)$, ${V_{2p}} = \mathbb{E}\left[ {\frac{1}{s}\left( {\sum_{i = 1}^N {\frac{{{{\hat d}_i}^2}}{{2{\beta _i}}}}  + \sum_{i = 1}^N {{\varpi _i}} } \right){\textbf{1}_{\left\{ {\sigma \left( t \right) = p} \right\}}}} \right]$, ${\hat d_i} = {d_i} - 4c - {o_i} - 1$, and ${{\textbf{1}_{\left\{ {\sigma \left( t \right) = p} \right\}}}}$ is the indicator function associated with the set $\left\{ {\sigma \left( t \right) = p} \right\}$.

We can prove $V>0$ by illustrating the positive definiteness of ${{V}_1}$, considering ${{V}_2}>0$. In fact, with some calculations, we have
\begin{align}
	{{V}_1} &= \sum\limits_{p = 1}^s {\mathbb{E}\left[ {\left( {{Z^T}{\Theta_p} Z} \right){\textbf{1}_{\left\{ {\sigma \left( t \right) = p} \right\}}}} \right]} \nonumber\\
	&= \mathbb{E}\left[ { {\sum\limits_{p = 1}^s {\Pi _P}{{Z^T}{\Theta_p} Z} } } \right]\nonumber\\
	& = \mathbb{E}\left[{Z^T}\left( {\begin{array}{*{20}{c}}
			{(\sum_{p=1}^s \Pi_p \mathcal{L}_p) \otimes P}&0\\
			0&{{I_N} \otimes Q}
	\end{array}} \right)Z \right]\nonumber\\
	&\ge \mathbb{E}\left[\sum_{p=1}^s\sum\limits_{i = 1}^N {\Pi_p \mathcal{L}_p \delta _i^T P\delta_i} + \sum\limits_{i = 1}^N {q_i^T} Qq_i \right] > 0,
\end{align}
 where the last inequality holds base on property $\mathop {\min }\limits_{x \ne 0,{1^T}x = 0} \frac{{{x^T}{\cal L}x}}{{{x^T}x}} = {\lambda _2}\left( {\cal L} \right)$ in \cite{olfati2004}.

 Computing the time derivative of ${{V_{1}}}$ in conjunction with (\ref{eq14}) and (\ref{eq15}) yields
\begin{align}\label{eq27} 
	&{{\dot V}_1} =\mathbb{E}\left[ {\sum\limits_{p = 1}^s {{\Pi _P}\left( {2{\delta ^T}\left( {{\cal L}\left( p \right) \otimes PA} \right)\delta } \right.} } \right.\nonumber\\
	&- c{\delta ^T}\left( {{{\cal L}^2}\left( p \right) \otimes \Gamma } \right)\delta   + c{\delta ^T}\left( {{{\cal L}^2}\left( p \right) \otimes \Gamma } \right)\delta \nonumber\\
	& + 2{\delta ^T}\left( {{\cal L}\left( p \right)Md \otimes \Gamma } \right)\tilde \xi+ 2{q^T}\left( {\frac{1}{s}{I_N} \otimes Q\left( {A - GC} \right)} \right)q\nonumber\\
	&\left. {\left. { + 2{q^T}\left( {\frac{1}{s}F \otimes QG} \right)\varepsilon } \right)} \right].
\end{align}

For the summation of the first two terms in (\ref{eq27}), we can obtain
\begin{align}\label{proof_eq_1}
	&\mathbb{E}\left[ {\sum\limits_{p = 1}^s {{\Pi _P}\left( {{\delta ^T}\left( {2{\cal L}\left( p \right) \otimes PA - c{{\cal L}^2}\left( p \right) \otimes \Gamma } \right)} \right.} \delta } \right]\nonumber\\
	&= \mathbb{E}\left[ {\sum\limits_{p = 1}^s {{\Pi _P}\left( {{\delta ^T}\left( {{\cal L}\left( p \right) \otimes \left( {{A^T}P + PA} \right)} \right)} \right.} \delta } \right.\nonumber\\
	&\quad \quad \quad \quad \left. { - \sum\limits_{p = 1}^s {{\Pi _P}\left( {{\delta ^T}\left( {c{{\cal L}^2}\left( p \right) \otimes \Gamma } \right)} \right.} \delta } \right]\nonumber\\
	&\le \mathbb{E}\bigg[ {{\delta ^T}( {\bar \Pi \hat {\cal L} \otimes( {{A^T}P + PA} )} )} \delta \nonumber\\
	&\quad \quad \quad \quad   - {\delta ^T}( {\frac{c}{s}{{\bar \Pi }}{( \sum\limits_{p = 1}^s {{\cal L}( p )}  )^2} \otimes \Gamma } )\delta  \bigg],
\end{align}
where the inequality holds following from the fact that ${\left[ {\sum_{p = 1}^s {{\cal L}\left( p \right)} } \right]^2} \le s\sum_{p = 1}^s {{{\cal L}^2}\left( p \right)}$.

Let $\xi  \triangleq  \left[  {\xi _1^T,...,\xi _N^T} \right]^T =  - \left( {{\cal L}\left( p \right) \otimes {I_n}} \right)\left( {q + \delta } \right)$. By using Young’s inequality Lemma, we have the following two inequalities
\begin{align}\label{yang1}
	&c{\delta ^T}\left( {{{\cal L}^2}\left( p \right)} \right.\left. { \otimes \Gamma } \right)\delta \nonumber\\
	&= c\left[ {{\xi ^T}\left( {{I_N} \otimes \Gamma } \right)\xi } \right. + 2{\xi ^T}\left( {{\cal L}\left( p \right) \otimes \Gamma } \right)q \left. { + {q^T}\left( {{{\cal L}^2}\left( p \right) \otimes \Gamma } \right)q} \right]\nonumber\\
	&\le 2c\left[ {{\xi ^T}\left( {{I_N} \otimes \Gamma } \right)\xi  + {q^T}\left( {{{\cal L}^2}\left( p \right) \otimes \Gamma } \right)q} \right],
\end{align}
\begin{align} \label{yang2}
	&2{\delta ^T}\left( {{\cal L}\left( p \right)Md \otimes \Gamma } \right){\tilde \xi }\nonumber\\
	&=  - 2{\xi ^T}\left( {d \otimes \Gamma } \right){\tilde \xi } - 2{q^T}\left( {{\cal L}\left( p \right)d \otimes \Gamma } \right){\tilde \xi }\nonumber\\
	&\le  - 2{\xi ^T}\left( {d \otimes \Gamma } \right){\tilde \xi } + {{\tilde \xi } ^T}\left( {{I_N} \otimes \Gamma } \right){\tilde \xi } \nonumber\\
	&+ {q^T}\left( {{{\cal L}^2}\left( p \right){{\bar d}^2} \otimes \Gamma } \right)q,
\end{align}
where the first equality holds because ${\cal L}\left( p \right) = {{\cal L}^T}\left( p \right)$, and $ - \left( {{\cal L}\left( p \right) \otimes {I_n}} \right)\delta  = \left( {{\cal L}\left( p \right) \otimes {I_n}} \right)q + \xi$.

%{\color{red}Following from the fact that ${\left[ {\sum\limits_{p = 1}^s {{\cal L}\left( p \right)} } \right]^2} \le s\sum\limits_{p = 1}^s {{{\cal L}^2}\left( p \right)}$, we get $ - \sum\limits_{p = 1}^s {{{\cal L}^2}\left( p \right)}  \le  - \frac{1}{{\rm{s}}}{\left[ {\sum\limits_{p = 1}^s {{\cal L}\left( p \right)} } \right]^2}$.}

Define ${e_{{\xi _i}}} \triangleq {\tilde{\xi} _i} - {\xi _i}$ with ${\xi _i} = \mathop \sum_{j = 1}^N {a_{ij}}\left( t \right)\left( {{{\hat x}_j} - {{\hat x}_i}} \right)$, and rewrite in a vector form  ${e_{{\xi }}}=\left[  {{e_{{\xi _1}}}^T,...,{e_{{\xi _N}}}^T} \right]^T$.
Substituting (\ref{yang1}) and (\ref{yang2}) into the summation of the first four terms of (\ref{eq27}) and considering (\ref{proof_eq_1}) yield
\begin{align}\label{v11}
	&\mathbb{E}\bigg[ {\sum\limits_{p = 1}^s {{\Pi _P}( {2{\delta ^T}( {{\cal L}(p) \otimes PA} )\delta } } }  - c{\delta ^T}( {{{\cal L}^2}(p) \otimes \Gamma } )\delta \nonumber\\
	&+ c{\delta ^T}( {{{\cal L}^2}(p)}{ \otimes \Gamma } )\delta  + {{2{\delta ^T}( {{\cal L}(p)Md \otimes \Gamma } )\tilde \xi } )} \bigg]\nonumber\\
	&\le \mathbb{E}\bigg[ {\sum\limits_{p = 1}^s {{\Pi _P}( {{\delta ^T}( {2{\cal L}(p) \otimes PA - c{{\cal L}^2}(p) \otimes \Gamma })}} }\delta \nonumber\\
	&+ 2c( {{\xi ^T}( {{I_N} \otimes \Gamma } )\xi  + {q^T}( {{{\cal L}^2}(p) \otimes \Gamma })q}) + {{\tilde \xi }^T}( {{I_N} \otimes \Gamma } )\tilde \xi \nonumber\\
	&- 2{\xi ^T}( {d \otimes \Gamma })\tilde \xi { { + {q^T}( {{{\cal L}^2}( p ){{\bar d}^2} \otimes \Gamma } )q} )} \bigg]\nonumber\\
	&\le \mathbb{E}\left[ {{\delta ^T}\left( {\bar \Pi \hat {\cal L} \otimes \left( {{A^T}P + PA} \right) - \frac{c}{s}\bar \Pi \hat {\cal L}{\lambda _2}\left( {\hat {\cal L}} \right) \otimes \Gamma } \right)} \right.\delta \nonumber\\
	&+ 2c{\xi ^T}\left( {{I_N} \otimes \Gamma } \right)\xi  - 2{\xi ^T}\left( {d \otimes \Gamma } \right)\tilde \xi \nonumber\\
	&+ s\left( {2c + {{\bar d}^2}} \right)\mathord{\buildrel{\lower3pt\hbox{$\scriptscriptstyle\smile$}} 
		\over \Pi } \lambda _M^2\left( {\hat {\cal L}} \right){q^T}\left( {{I_N} \otimes \Gamma } \right)q\left. { + {{\tilde \xi }^T}\left( {{I_N} \otimes \Gamma } \right)\tilde \xi } \right]\nonumber\\
	&= \mathbb{E}\left[ {{\delta ^T}\left( {\bar \Pi \hat {\cal L} \otimes \left( {{A^T}P + PA} \right) - \frac{c}{s}\bar \Pi \hat {\cal L}{\lambda _2}\left( {\hat {\cal L}} \right) \otimes \Gamma } \right)} \right.\delta \nonumber\\
	&+ 2c{\xi ^T}\left( {{I_N} \otimes \Gamma } \right)\xi  + s\left( {2c + {{\bar d}^2}} \right)\mathord{\buildrel{\lower3pt\hbox{$\scriptscriptstyle\smile$}} 
		\over \Pi } \lambda _M^2\left( {\hat {\cal L}} \right){q^T}\left( {{I_N} \otimes \Gamma } \right)q\nonumber\\
	&- 2{{\tilde \xi }^T}\left( {d \otimes \Gamma } \right)\tilde \xi  + 2e_\xi ^T\left( {d \otimes \Gamma } \right)\tilde \xi \left. { + {{\tilde \xi }^T}\left( {{I_N} \otimes \Gamma } \right)\xi } \right],
\end{align}
where ${\cal L}\left( p \right)M = {\cal L}\left( p \right) = M{\cal L}\left( p \right)$ is used in the first inequality.

Consider the following formula
\begin{align}\label{v11v}
	&2c{\xi ^T}\left( {{I_N} \otimes \Gamma } \right)\xi  - 2{{\tilde \xi }^T}\left( {d \otimes \Gamma } \right)\tilde \xi  + 2e_\xi ^T\left( {d \otimes \Gamma } \right)\tilde \xi \nonumber\\
	&\quad + {{\tilde \xi }^T}\left( {{I_N} \otimes \Gamma } \right)\tilde \xi \nonumber\\
	&\le 2c{\xi ^T}\left( {{I_N} \otimes \Gamma } \right)\xi  - {{\tilde \xi }^T}\left( {d \otimes \Gamma } \right)\tilde \xi  + e_\xi ^T\left( {d \otimes \Gamma } \right){e_\xi }\nonumber\\
	&\quad + {{\tilde \xi }^T}\left( {{I_N} \otimes \Gamma } \right)\tilde \xi \nonumber\\
	&= 2c{{\tilde \xi }^T}\left( {{I_N} \otimes \Gamma } \right)\tilde \xi  + 2ce_\xi ^T\left( {{I_N} \otimes \Gamma } \right){e_\xi } - 4c{{\tilde \xi }^T}\left( {{I_N} \otimes \Gamma } \right){e_\xi }\nonumber\\
	&\quad - {{\tilde \xi }^T}\left( {d \otimes \Gamma } \right)\tilde \xi  + e_\xi ^T\left( {d \otimes \Gamma } \right){e_\xi } + {{\tilde \xi }^T}\left( {{I_N} \otimes \Gamma } \right)\tilde \xi \nonumber\\
	&\le e_\xi ^T\left( {\left( {d + 4c{I_N}} \right) \otimes \Gamma } \right){e_\xi } \nonumber\\
	&\quad - {{\tilde \xi }^T}\left( {\left( {d - \left( {4c + 1} \right){I_N}} \right) \otimes \Gamma } \right)\tilde \xi, 
\end{align}
in which the first and second inequalities hold because of the Young’s Inequality Lemma; the equality is satisfied due to the fact that ${\xi} = {\tilde{\xi} } - {e_{{\xi }}}$. 

Meanwhile, we have ${e_{\xi i}} = \sum\limits_{j = 1,j \ne i}^N {{a_{ij}}\left( {\sigma \left( t \right)} \right)} \left( {{m_j} - {m_i}} \right)$ for $t=t_k^i$, with the definition of ${e_{{\xi _i}}}$ and some mathematical calculations. Then, we can obtain
\begin{align}\label{xin1}
	&\mathbb{E}\left[ {\sum\limits_{p = 1}^s {{\Pi _p}\left\| {e_{{\xi _i}}^TPB} \right\|} } \right]\nonumber\\
	&= \mathbb{E}\left[ {\sum\limits_{p = 1}^s {{\Pi _p}\left\| {\sum\limits_{j = 1,j \ne i}^N {{a_{ij}}} \left( {\sigma \left( t \right)} \right){{\left( {{m_j} - {m_i}} \right)}^T}PB} \right\|} } \right]\\
	&\le \mathbb{E}\left[ {\mathord{\buildrel{\lower3pt\hbox{$\scriptscriptstyle\smile$}} 
			\over \Pi } \left( {{{\hat {\cal L}}_{ii}}\left\| {m_i^TPB} \right\| + \sum\limits_{j = 1,j \ne i}^N {{{\hat a}_{ij}}} \left\| {m_j^TPB} \right\|} \right)} \right]\nonumber\\	&\le \mathbb{E}\left[ {\mathord{\buildrel{\lower3pt\hbox{$\scriptscriptstyle\smile$}} 
			\over \Pi } N\left( {\hat {\cal L}_{ii}^2{{\left\| {m_i^TPB} \right\|}^2} + \sum\limits_{j = 1,j \ne i}^N {{{\hat a}_{ij}}} {{\left\| {m_j^TPB} \right\|}^2}} \right)} \right],\nonumber
\end{align}
where the last inequality is satisfied because of the fact that ${\bigg( {\sum\limits_{i = 1}^N {{x_i}} } \bigg)^2} \le N\sum\limits_{i = 1}^N {x_i^2} $. Therefore, for the first term in (\ref{v11v}), we have
\begin{align}
	&\mathbb{E}\left[\sum\limits_{p = 1}^s{ {{\Pi _p}\left( {e_\xi ^T\left( {\left( {d + 4c{I_N}} \right) \otimes \Gamma } \right){e_\xi }} \right)}} \right] \nonumber\\
	&=\mathbb{E}\left[\sum\limits_{p = 1}^s{ {{\Pi _p}\left( {\sum\limits_{i = 1}^N {\left( {{d_i} + 4c} \right){{\left\| {e_{{\xi _i}}^TPB} \right\|}^2}} } \right)}} \right] \nonumber\\
	&\le \mathbb{E}\bigg[ {\mathord{\buildrel{\lower3pt\hbox{$\scriptscriptstyle\smile$}} 
			\over \Pi } N\sum\limits_{i = 1}^N {( {{{\bar d}_i} + 4c} )} ( {\hat {\cal L}_{ii}^2{{\| {m_i^TPB} \|}^2}} } \nonumber\\
	&\quad \quad  { { + \sum\limits_{j = 1,j \ne i}^N {\hat a_{ij}} {{\| {m_j^TPB} \|}^2}} )} \bigg]. \nonumber
\end{align}

Further computations reveal that
\begin{align}\label{xin2}
&\mathbb{E}\left[\sum\limits_{p = 1}^s{ {{\Pi _p}\left( {e_\xi ^T\left( {\left( {d + 4c{I_N}} \right) \otimes \Gamma } \right){e_\xi }} \right)}} \right] \nonumber\\
&\le \mathbb{E}\bigg[ {\mathord{\buildrel{\lower3pt\hbox{$\scriptscriptstyle\smile$}} 
		\over \Pi } N\sum\limits_{i = 1}^N {( {( {{{\bar d}_i} + 4c} )\hat {\cal L}_{ii}^2 + \sum\limits_{j = 1,j \ne i}^N {( {{{\bar d}_j} + 4c} )\hat a_{ji}} } )} } \nonumber\\
&\quad \quad  { \times {{\| {m_i^TPB} \|}^2}} \bigg]\nonumber\\
&\le \mathbb{E}\left[ {\mathord{\buildrel{\lower3pt\hbox{$\scriptscriptstyle\smile$}} 
		\over \Pi } N{\rm{\tilde c}}\sum\limits_{i = 1}^N {\left( {\hat {\cal L}_{ii}^2 + \sum\limits_{j = 1,j \ne i}^N {\hat a_{ji}^2} } \right)} {{\left\| {m_i^TPB} \right\|}^2}} \right]\nonumber\\
&\le \mathbb{E}\left[ {{\rm{\tilde c}}\mathord{\buildrel{\lower3pt\hbox{$\scriptscriptstyle\smile$}} 
		\over \Pi } N\sum\limits_{i = 1}^N {\left( {{N^2} + N} \right)} {{\left\| {m_i^TPB} \right\|}^2}} \right],
\end{align}
where the second inequality holds by using the definition that ${\rm{\tilde c}} = \mathop {\max }\limits_{i = 1,2,...,N} \left\{ {{{\bar d}_i} + 4c} \right\}$.

The last two terms of (\ref{eq27}) satisfy that
\begin{align}\label{v12}
	&\mathbb{E}\left[ {2{q^T}\left[ {{I_N} \otimes \left( {QA - QGC} \right)} \right]q} \right. \left. { + 2{q^T}\left( {F \otimes QG} \right)\varepsilon } \right]\nonumber\\
	&\le \mathbb{E}\left[ {{q^T}\left( {{I_N} \otimes \left( {{A^T}Q + QA - QGC - {{\left( {QGC} \right)}^T}} \right)} \right)} \right.q\nonumber\\
	&\quad \left. { + {q^T}\left( {F \otimes QG} \right){{\left( {F \otimes QG} \right)}^T}q + {\varepsilon ^T}\varepsilon } \right]\nonumber\\
	&\le \mathbb{E}\bigg[ {{q^T}( {{I_N} \otimes ( {{A^T}Q + QA - QGC - {{( {QGC} )}^T}} } } \nonumber\\
	&\quad  { { { + {\lambda _M}( {F{F^T}} )QG{G^T}{Q^T}} )} )q + {\varepsilon ^T}\varepsilon } \bigg],
\end{align}
where the first inequality is satisfied because of the Young’s Inequality Lemma.

Calculating the time derivative of ${V_2}$ yields
\begin{align}\label{v2}
	{{\dot V}_2}
	&= \mathbb{E}\bigg[ {\sum\limits_{i = 1}^N {\frac{1}{{{\beta _i}}}} {{\hat d}_i}{{\dot d}_i}}  + \sum\limits_{i = 1}^N {\bigg( { - {\eta _i}{\varpi _i}} } \nonumber\\
	&\quad \qquad  { { + {\varsigma _i}\left[ {{o_i}{\upsilon _i}{\tilde \xi } _i^T\Gamma {{\tilde \xi } _i} - m_i^T\Gamma {m_i}} \right]} \bigg)} \bigg].
\end{align}

Substituting (\ref{v11}), (\ref{v11v}), (\ref{xin2}) and (\ref{v12}) into (\ref{eq27}) and considering (\ref{v2}), we have the derivative of the candidate Lyapunov function yield 
\begin{align}\label{vo1}
	\dot V &\le \mathbb{E}\bigg[ {\delta ^T}( \bar \Pi \hat {\cal L} \otimes ( {{A^T}P + PA} ) - \frac{c}{s}\bar \Pi \hat {\cal L}{\lambda _2}( {\hat {\cal L}} ) \otimes \Gamma  ) \delta \nonumber\\
	&+ {m^T}( \rho {I_N} \otimes \Gamma )m - \sum\limits_{i = 1}^N {\varsigma _i}m_i^T\Gamma {m_i}  + {q^T}( {{I_N} \otimes \Psi } )q  \nonumber\\
	&+ {\varepsilon ^T}\varepsilon  - {{\tilde \xi }^T}( ( {d - ( {4c + 1 + o} ){I_N}} ) \otimes \Gamma  )\tilde \xi  + \sum\limits_{i = 1}^N {\frac{1}{{{\beta _i}}}} {{\hat d}_i}{{\dot d}_i}\nonumber\\
	& { - {{\tilde \xi }^T}( {o{I_N} \otimes \Gamma } )\tilde \xi  + \sum\limits_{i = 1}^N {{\varsigma _i}{o_i}{\upsilon _i}\tilde \xi _i^T\Gamma {{\tilde \xi }_i}}  - \sum\limits_{i = 1}^N {{\eta _i}{\varpi _i}} } \bigg].
\end{align}

Utilizing the update law provided in (\ref{eq11}), we proceed with the subsequent analysis as follows

(1) if ${{d}_i} < {{\bar d}_i}$, it follows from (\ref{eq11}) that ${{\dot d}_i} = {{\beta _i}}{{\tilde \xi } _i}^T\Gamma {{\tilde \xi } _i}$. Then, we can get 
\begin{equation}\label{if1}
	\sum\limits_{i = 1}^N {\frac{1}{{{\beta _i}}}} {{\hat d}_i}{{\dot d}_i} - {{\tilde \xi } ^T}\left( {\hat d \otimes \Gamma } \right){\tilde \xi } = 0, \text{if}\ {{d}_i} < {{\bar d}_i}.
\end{equation}

(2) if ${{d}_i} \ge {{\bar d}_i}$, we can obtain from (\ref{eq11}) that ${{\dot d}_i} = 0$. Now, we have
\begin{equation}\label{if2}
	\sum\limits_{i = 1}^N {\frac{1}{{{\beta _i}}}} {{\hat d}_i}{{\dot d}_i} - {{\tilde \xi } ^T}\left( {\hat d \otimes \Gamma } \right){\tilde \xi } < 0,  \text{if}\ {{d}_i} \ge {{\bar d}_i}.
\end{equation}
Combining (\ref{if1}) and (\ref{if2}) yields
\textbf{\begin{equation}\label{if22}
		\sum\limits_{i = 1}^N {\frac{1}{{{\beta _i}}}} {{\hat d}_i}{{\dot d}_i} - {{\tilde \xi } ^T}\left( {\hat d \otimes \Gamma } \right){\tilde \xi } \leq 0.
\end{equation}}

According to ${\upsilon _i} \geq \frac{1}{{{\rho}}}$ and ${\hat d_i} = {d_i} - 4c - 1 - {o_i}$, and substituting (\ref{if22}) into (\ref{vo1}) yield
\begin{align} \label{vo2}
	\dot V &\le \mathbb{E}\bigg[ {\delta ^T}( \bar \Pi \hat {\cal L} \otimes ( {A^T}P + PA ) - \frac{c}{s}\bar \Pi {\lambda _2}( \hat {\cal L} )\hat {\cal L} \nonumber\\
	&  \otimes \Gamma  )\delta  + {q^T}( {I_N} \otimes \Psi  )q + \sum\limits_{i = 1}^N ( \rho  - {\varsigma _i} )m_i^T\Gamma {m_i} \nonumber\\
	&+ {\varepsilon ^T}\varepsilon   - \sum\limits_{i = 1}^N ( \rho  - {\varsigma _i} ){o_i}{\upsilon _i}\tilde \xi _i^T\Gamma {\tilde \xi }_i  - \sum\limits_{i = 1}^N {\eta _i}{\varpi _i}  \bigg].
\end{align}

Note that during the time period $t\in [t_k^i, t_{k+1}^i)$, the event-triggered condition given in (\ref{eq9}) is not satisfied. With ${\rho _i} - {\varsigma _i} > 0$ in (\ref{eq17}), we can obtain
\begin{align}
	&\sum\limits_{i = 1}^N {\left( {{\rho _i} - {\varsigma _i}} \right)} m_i^T\Gamma {m_i} - \sum\limits_{i = 1}^N {\left( {{\rho _i} - {\varsigma _i}} \right)} {o_i}{\upsilon _i}{{\tilde \xi } _i}^T\Gamma {{\tilde \xi } _i} \nonumber\\
	&\qquad  < \sum\limits_{i = 1}^N {\frac{{{\rho _i} - {\varsigma _i}}}{{{\iota _i}}}{\varpi _i}}.
\end{align}

According to ${\eta _i} - \frac{{{\rho} - {\varsigma _i}}}{{{\iota _i}}}>0$ in (\ref{eq17}), it follows from (\ref{vo2}) that
\begin{align}
	\dot V &\le \mathbb{E}\bigg[ {{\delta ^T}( {\hat {\cal L} \otimes \bar \Pi ( {( {{A^T}P + PA} ) - \frac{c}{s}{\lambda _2}( {\hat {\cal L}} )\Gamma } )} )} \delta \nonumber\\
	&+ {q^T}( {{I_N} \otimes \Psi } )q + \tau  {  - \sum\limits_{i = 1}^N {( {{\eta _i} - \frac{{\rho  - {\varsigma _i}}}{{{\iota _i}}}} )} {\varpi _i} } \bigg]\nonumber\\
	&\le \mathbb{E}\left[ {{Z^T}\left( {\left( {\begin{array}{*{20}{c}}
					{\hat {\cal L}}&0\\
					0&{{I_N}}
			\end{array}} \right) \otimes \Lambda } \right)Z} \right] + \tau \nonumber\\
	&\le  - \kappa \mathbb{E}\left[ {{Z^T}\left( {\left( {\begin{array}{*{20}{c}}
					{\hat {\cal L}}&0\\
					0&{{I_N}}
			\end{array}} \right) \otimes \tilde \Theta } \right)Z} \right] + \tau \nonumber\\
	&\le  - \kappa {V_1} + \tau,
\end{align}
where the third inequality holds because of $\Lambda  + \kappa \tilde \Theta  < 0$ in (\ref{eq17}).

Since $V = {V_1} + {V_2}$, according to the comparison Lemma, we have
\begin{equation}\label{proof_ineq_2}
	{{\dot V}_1} \le  - \kappa {{V_1}} + \tau .
\end{equation}

Integrating both sides of (\ref{proof_ineq_2}) from $s=t_0$ to $s=t$ yields 
\begin{equation}
	{V_1} \le {e^{ - \kappa \left( {t - {t_0}} \right)}}{{V_1}\left( {{t_0}} \right)} + \tau \int_{{t_0}}^t {{e^{ - \kappa (t - s)}}} ds.
\end{equation}

From (\ref{eq18}), a deduction of $\chi I < \left( {\begin{array}{*{20}{c}}
		{\hat {\cal L} \otimes P}&0\\
		0&{{I_N} \otimes Q}
\end{array}} \right)$ is achieved. Then, we can obtain

\begin{align}
	\mathbb{E}\left[ {{{\left\| {Z} \right\|}^2}} \right] &\le\frac{1}{\chi }\mathbb{E}\left[ {{V_1}} \right]\nonumber\\
	&\le \frac{1}{\chi }{e^{ - \kappa \left( {t - {t_0}} \right)}}\mathbb{E}\left[ {{V_1}\left( {{t_0}} \right)} \right] + \frac{\tau }{{\chi \kappa }} - \frac{\tau }{{\chi \kappa }}{e^{ - \kappa \left( {t - {t_0}} \right)}}.
\end{align}

As $t \to \infty $, we have
\begin{align}\label{eq42}
	\mathbb{E}\left[ {\left\| {\delta } \right\|} \right] \le E\left[ {\left\| {Z} \right\|} \right] \le \sqrt {\frac{\tau }{{\chi \kappa }}}, \nonumber\\
	\mathbb{E}\left[ {\left\| {q} \right\|} \right] \le \mathbb{E}\left[ {\left\| {Z} \right\|} \right] \le \sqrt {\frac{\tau }{{\chi \kappa }}} .
\end{align}

It follows from (\ref{eq42}) that the consensus error ${\delta}$ converges to the set ${\bar \ell }$ as $t \to \infty $ in the mean square sense. The proof is completed.

	\begin{remark}
To overcome the complexities arising from deception attacks and the intricate interplay between Markovian randomly switching topologies and agents' dynamics, this study proposes an observer-based distributed secure control protocol, ensuring the system's attainment of both stability and security.
\end{remark}
%\end{IEEEproof}
\begin{theorem}
	The Zeno behavior of the MASs (\ref{eq2}) can be effectively eliminated under the observer-based distributed secure controller with control gain matrix given in Theorem \ref{thm1}, ensuring a minimum time interval between consecutive events.
\end{theorem}
%\begin{proof}

\textit{Proof}: Supposing that ${t_k^i}$ and ${t_{k + 1}^i}$ are any two consecutive triggering moments on interval $\left[ {{t_0},t} \right)$. For $t \in \left( {t_k^i,t_{k + 1}^i} \right)$, based on the dynamic event-triggered condition (\ref{eq9}) and the proof of Lemma $1$, we have
\begin{equation}
	\left( {m_i^T\Gamma {m_i} - {o_i}{\upsilon _i}{\tilde \xi } _i^T\Gamma {{\tilde \xi } _i}} \right) < \frac{{{\varpi _i}\left( 0 \right)}}{{{\iota _i}}}{e^{ - \left( {{\eta _i} + \frac{{{\varsigma _i}}}{{{\iota _i}}}} \right)t}}.
\end{equation}

Taking derivative of ${m_i} = {\hat x_i}(t_k^i) - {\hat x_i}$ yields
\begin{align}\label{eq44}
	{{\dot m}_i} =  - {\dot {\hat x}}_{i} &  =  - \left( {A{{\hat x}_i} + {{d}_i}BK{{\tilde \xi } _i}} \right.\left. { + G\left( {{{\mathord{\buildrel{\lower3pt\hbox{$\scriptscriptstyle\smile$}} 
						\over y} }_i} - {{\hat y}_i}} \right)} \right)\nonumber\\
	&  = A{m_i} - A{{\hat x}_i}\left( {t_k^i} \right) - {{d}_i}BK{{\tilde \xi } _i}+ GC{q_i} - G{\alpha}{\varepsilon _i}\nonumber\\
	&\buildrel \Delta \over = A{m_i} + {D_i}\left( {t_k^i} \right) + {\Omega _i},
\end{align}
where ${D_i}\left( {t_k^i} \right) =  - A{\hat x_i}\left( {t_k^i} \right) - {d_i}BK{\xi _i}\left( {t_k^i} \right), {\Omega _i} = GC{q_i} - G{\alpha}{\varepsilon _i}$. It follows from (\ref{eq44}) that
\begin{align}
	&\int_{t_k^i}^{\rm{t}} {{e^{ - As}}\left( {{{\dot m}_i} - A{m_i}} \right)} ds= \int_{t_k^i}^{\rm{t}} {{e^{ - As}}\left( {{D_i}\left( {t_k^i} \right) + {\Omega _i}} \right)} ds.
\end{align}

Further, it is equivalent to
\begin{equation}
	{e^{ - At}}{m_i} - {e^{ - At_k^i}}{m_i}\left( {t_k^i} \right) = \int_{t_k^i}^{\rm{t}} {{e^{ - As}}\left( {{D_i}\left( {t_k^i} \right) + {\Omega _i}} \right)} ds.
\end{equation}

Considering the fact ${m_i}\left( {t_k^i} \right) = 0$, we have
\begin{equation}
	{m_i} = \int_{t_k^i}^t {{e^{A\left( {t - s} \right)}}\left( {{D_i}\left( {t_k^i} \right) + {\Omega _i}} \right)ds}.
\end{equation}

According to Theorem \ref{thm1}, it has been shown that there exist positive constants ${\theta _1}$, ${\theta _2}$, and $\nu$ such that the expected value of $q$ approaches ${\theta _1}$ with a convergence rate of $\nu$ within the interval $[t_k^i, t)$. Specifically, $\mathbb{E}\left[\parallel {q} \parallel\right] \le {\theta _2}{e^{ - \nu \left( {t - t_k^i} \right)}}\mathbb{E}\left[ \parallel{q\left( {t_k^i} \right)}\parallel \right] + {\theta _1}$. Furthermore, considering ${\Omega _i} = GC{q_i} - GF{\varepsilon _i}$, there exists a positive constant $\vartheta $ such that the following inequality holds
\begin{align}\label{etin1}
	&\mathbb{E}\left[ {\parallel {m_i}\parallel } \right] \nonumber\\
	&\le \mathop \smallint \nolimits_{t_k^i}^t {e^{\parallel A\parallel (t - s)}}\left( {\mathbb{E}\left[ {\parallel {D_i}\left( {t_k^i} \right)\parallel } \right] + {\theta _2}{e^{ - \nu \left( {t - t_k^i} \right)}}} \right.\nonumber\\
	&\times \left. {\parallel GC\parallel \left( {\mathbb{E}\left[ {\parallel q\left( {t_k^i} \right)\parallel } \right] + {\theta _1}} \right) + \frac{\tau }{N}\parallel G\parallel } \right)ds\nonumber\\
	&  \le {\rm{\Phi }}\left( {t_k^i} \right)\mathop \smallint \nolimits_{t_k^i}^t {e^{\parallel A\parallel (t - s)}}ds - \theta \left( {t_k^i} \right),
\end{align}
where ${\rm{\Phi }}\left( {t_k^i} \right) = \mathbb{E}\left[ {\parallel {D_i}\left( {t_k^i} \right)\parallel } \right] + {\theta _2}\parallel GC\parallel \left( {\mathbb{E}\left[ {\parallel q\left( {t_k^i} \right)\parallel } \right] + {\theta _1}} \right) + \frac{\tau }{N}\parallel G\parallel  + \vartheta$, $\theta \left( {t_k^i} \right) > 0$.

1) If $\left\| A \right\| \ne 0$,  we can obtain from (\ref{etin1}) that
\begin{align}\label{eq49}
	\mathbb{E}[ {\parallel m_i^TPB\parallel } ] &  \le \bigg( {\frac{{{\rm{\Phi }}( {t_k^i} )}}{{\| A \|}}( {{e^{\| A \|( {t - t_k^i} )}} - 1} )} \nonumber\\
	& { - \theta ( {t_k^i} )}  \bigg)\| {PB} \|.
\end{align}

Using the dynamic event-triggered strategy (\ref{eq9}), the subsequent triggering instant ${t_{k + 1}^i}$ is determined by
\begin{align}\label{eq50}
	{\left\| {m_i^T\left( {t_{k + 1}^i} \right)PB} \right\|^2} &\ge {o_i}{\upsilon _i}{\left\| {\xi _i^T\left( {t_{k + 1}^i} \right)PB} \right\|^2}\nonumber\\
	&+ \frac{{{\varpi _i}\left( 0 \right)}}{{{\iota _i}}}{e^{ - \left( {{\eta _i} + \frac{{{\varsigma _i}}}{{{\iota _i}}}} \right)t_{k + 1}^i}}.
\end{align}

Then, we have
\begin{align}\label{eq51}
	0 &< \mathbb{E}\left[ {\frac{{{\varpi _i}\left( 0 \right)}}{{{\iota _i}}}{e^{ - \left( {{\eta _i} + \frac{{{\varsigma _i}}}{{{\iota _i}}}} \right)t_{k + 1}^i}}} \right] \nonumber\\
	&\le \mathbb{E}\left[ {{o_i}{\upsilon _i}{{\left\| {\xi _i^T\left( {t_{k + 1}^i} \right)PB} \right\|}^2} + \frac{{{\varpi _i}\left( 0 \right)}}{{{\iota _i}}}{e^{ - \left( {{\eta _i} + \frac{{{\varsigma _i}}}{{{\iota _i}}}} \right)t_{k + 1}^i}}} \right] \nonumber\\
	&\le \mathbb{E}\left[ {{{\left\| {m_i^T\left( {t_{k + 1}^i} \right)PB} \right\|}^2}} \right] \nonumber\\
	&\le \mathbb{E}\bigg[ {{\bigg( {\frac{{{\rm{\Phi }}( {t_{k + 1}^i} )}}{{\| A \|}}( {{e^{  \| A \|( {t_{k + 1}^i - t_k^i} )}} - 1} ) - \theta ( {t_k^i} )} \bigg)^2}} \nonumber\\
	&{ \times {{\| {PB} \|}^2}} \bigg],
\end{align}
where the third and forth inequalities hold because of (\ref{eq50}) and (\ref{eq49}), respectively.

The time interval between two consecutive events satisfies
\begin{align} \label{jiange}
	&\Delta _k^i = t_{k + 1}^i - t_k^i \ge \frac{1}{{\| A \|}}\ln \bigg\{ {\frac{{\| A \|}}{{{\rm{\Phi }}( {t_{k + 1}^i} )}}} \nonumber\\
	&\times \bigg( {\frac{{\sqrt {{o_i}{\upsilon _i}{{\| {\xi _i^T( {t_{k + 1}^i} )PB} \|}^2} + \frac{{{\varpi _i}( 0 )}}{{{\iota _i}}}{e^{ - ( {{\eta _i} + \frac{{{\varsigma _i}}}{{{\iota _i}}}} )t_{k + 1}^i}}} }}{{\| {PB} \|}}} \nonumber\\
	&{ { + \theta ( {t_k^i} )} \bigg) + 1} \bigg\} > 0.
\end{align}

The proof presented demonstrates that $\Delta _k^i > 0$, regardless of the finite horizon considered.

To demonstrate this, we assume the contrary, i.e., $\sum\nolimits_{k = 0}^\infty  {\Delta _k^i}$ will be convergent, which mean $\mathop {\lim }\nolimits_{m \to \infty } \sum\nolimits_{k = 0}^m {\Delta _k^i} $ will be convergent. From (\ref{eq51}), we have

\begin{align}
	&\mathop {\lim }\limits_{k \to \infty } \mathbb{E}\Bigg[\frac{{{\varpi _i}\left( 0 \right)}}{{{\iota _i}}}{e^{ - \left( {{\eta _i} + \frac{{{\varsigma _i}}}{{{\iota _i}}}} \right)t_{k + 1}^i}}\Bigg] \nonumber\\
	&\le \mathop {\lim }\limits_{k \to \infty } \mathbb{E}\Bigg[{e^{\left( {{\eta _i} + \frac{{{\varsigma _i}}}{{{\iota _i}}}} \right)t_k^i}}\Bigg(\frac{{\Phi \left( {t_{k + 1}^i} \right)}}{{\left\| A \right\|}}\left( {{e^{\left\| A \right\|\left( {t_{k + 1}^i - t_k^i} \right)}} - 1} \right) \nonumber\\
	&- \theta \left( {t_k^i} \right){\Bigg)^2}{\left\| {PB} \right\|^2} \Bigg]
\end{align}
which implies that $\left[ {{\varpi _i}\left( 0 \right)/{\iota _i}} \right] \le 0$, which contradicts the fact that $\left[ {{\varpi _i}\left( 0 \right)/{\iota _i}} \right] > 0$. Thus, $\mathop {\lim }\nolimits_{k \to \infty } t_k^i = \infty $. Then, $\mathop {\lim }\nolimits_{m \to \infty } \sum\nolimits_{k = 0}^m {\Delta _k^i}  = \mathop {\lim }\nolimits_{k \to \infty } t_{k + 1}^i - {t_0} = \infty $, which is divergent, the assumption is invalid. This excludes the possibility of Zeno behavior.

2) If $\left\| A \right\| = 0$, we have
\begin{align}
	\mathbb{E}\left[ {\parallel {m_i}PB\parallel } \right] &  \le {\rm{\Phi }}\left( {t_k^i} \right)\mathop \smallint \nolimits_{t_k^i}^t ds - \theta \left( {t_k^i} \right).
\end{align}

Similar to 1), the conclusion is drawn. The proof is completed.
%\end{proof}
\begin{remark}
The distributed dynamic event-triggered mechanism, introduced in this study, integrates the benefits of both distributed and dynamic event-triggered mechanisms. This integration offers an efficient, flexible, and reliable approach to communication and collaboration among MASs. Furthermore, it is demonstrated that the proposed event-triggered mechanism is devoid of Zeno behavior.
\end{remark}
\begin{theorem}\label{thrm3}
	Under Assumptions \ref{assumption1}-\ref{assumption3}, given ${o_i} \in \left( {0,1} \right)$, ${\varpi _i}(0) >0 $, $\kappa  > 0$, $\chi  > 0$, considering the predefined maximum attack threshold of $\tau $ and the probability matrices $F$, the error system (\ref{eq15}) converges into the set ${\bar \ell }$ exponentially in the mean square sense, if there exist positive definite matrices $P$ and $Q$ which fulfill the subsequent inequalities
	\begin{equation}%\label{eq17}
		\bar \Pi \left( {\left( {{A^T}P + PA} \right) - \frac{c}{s}{\lambda _2}\left( {\hat {\cal L}} \right)\Gamma } \right) + \kappa P < 0
	\end{equation}
	%\scriptsize{ \setlength{\arraycolsep}{1.2pt}
		\begin{align}%\label{eq17}
			\left( {\begin{array}{*{20}{c}}
					{\bar \Psi }&{\lambda _M^{\frac{1}{2}}\left( {F{F^T}} \right)X}\\
					*&{ - {I_n}}
			\end{array}} \right) < 0
		\end{align}
		\begin{equation}%\label{eq17}
			{\chi {I_n} - {\lambda _2}\left( {\hat {\cal L}} \right)P} < 0
		\end{equation}
		\begin{equation}%\label{eq17}
			{\chi {I_n} - Q}  < 0
		\end{equation}
		where $\bar \Psi  = \left( {{A^T}Q + QA - XC - {{\left( {XC} \right)}^T}+ s\left( {2c + {{\bar d}^2}} \right)} \right.$ \\ $\left. { {\mathord{\buildrel{\lower3pt\hbox{$\scriptscriptstyle\smile$}} 
					\over \Pi }}\lambda _M^2\left( {\hat {\cal L}} \right)\Gamma + \kappa Q } \right)$.
		Further, the controller gain matrix is calculated with $K = {B^T}P$ and observer gain matrix is calculated with $G = {Q^{ - 1}}X$.
	\end{theorem}
	%\begin{proof}
	
	\textit{Proof}: Substituting the definition of matrices $\Lambda $ and $\tilde \Theta $ into  (\ref{eq17}) gives
	\begin{align}\label{eq58}
		\scriptsize{ \setlength{\arraycolsep}{1.2pt}
			\left( {\begin{array}{*{20}{c}}
					{\bar \Pi \left( {{A^T}P + PA - \frac{c}{s}{\lambda _2}\left( {\hat {\cal L}} \right)\Gamma } \right) + \kappa P}&0\\
					0&{\Psi  + \kappa Q}
			\end{array}} \right) < 0}.
	\end{align}
	
By using Schur's Complement Lemma, we can get (\ref{eq58}) is equivalent to
	\begin{equation}\label{eq59}
		\bar \Pi \left( { {{A^T}P + PA}  - \frac{c}{s}{\lambda _2}\left( {\hat {\cal L}} \right)\Gamma } \right) + \kappa P < 0
	\end{equation}
	and
	\begin{equation}\label{eq60}
		\Psi  + \kappa Q < 0.
	\end{equation}
	
	Therefore, the controller gain matrix is determined by $K = {B^T}P$, where $P$ represents a positive-definite solution to the linear matrix inequality given in (\ref{eq59}). Let $X = QG$. Substituting the definition of matrix ${{\Psi}}$ into (\ref{eq60}) and using the Schur's Complement Lemma, we can get
	\begin{align}\label{eq61}
		\left( {\begin{array}{*{20}{c}}
				{\bar \Psi }&{\lambda _M^{\frac{1}{2}}\left( {F{F^T}} \right)X}\\
				*&{ - {I_n}}
		\end{array}} \right) < 0.
	\end{align}
	
	Substituting the definition of matrix $\hat \Theta $ into (\ref{eq18}) yields
	\begin{equation}\label{eq62}
\left( {\begin{array}{*{20}{c}}
		{\chi {I_n} - {\lambda _2}\left( {\hat {\cal L}} \right)P}&0\\
		0&{\chi {I_n} - Q}
\end{array}} \right) < 0.
	\end{equation}
	
	Again using the Schur's Complement Lemma, (\ref{eq62}) is equivalent to
	\begin{equation}\label{eq63}
		{\chi {I_n} - {\lambda _2}\left( {\hat {\cal L}} \right)P} <0,
			\end{equation}
	and
	 \begin{equation}\label{eq64}
		{{\chi {I_n} - Q}}< 0.
	\end{equation}
	
	By combining the preceding proof with the established results, it can be concluded that if equations (\ref{eq59}), (\ref{eq61}), (\ref{eq63}), and (\ref{eq64}) are satisfied, they are equivalent to the fulfillment of equations (\ref{eq17}) and (\ref{eq18}). Furthermore, the observer gain matrix can be determined as $G = Q^{-1}X$. The proof is now complete.
	%	
	%\end{proof}
	\begin{remark}
The design of controller and observer in this paper embraces a distributed paradigm, leveraging solely local information for the computation of the gain matrix within the MASs. This approach facilitates the design and implementation process of the entire system, enhancing its simplicity and accessibility.
	\end{remark}
\begin{table}[h]
	\caption{PARAMETER VALUES IN SIMULATION}
	\label{table}
	\setlength{\tabcolsep}{2pt}
	\begin{tabular}{|p{150pt}|p{93pt}|}
		\hline
		Symbol& 
		value \\
		\hline
		The upper bound of deception attack energy $\tau $ &
		$0.02$ \\
		The distribution $\Pi $ of the switching process&
		$\left[ {\begin{array}{*{20}{c}}
				{2/3}&{1/3}
		\end{array}} \right]$ \\
		Deception attack probability&
		$\begin{array}{l}
			\text{diag}\left\{ {0.32,0.24,0.30,} \right.\\
			0.42,0.27,0.32,0.25,\\
			\left. {0.23,0.39,0.28} \right\}
		\end{array}$ \\
		$c$ & 
		$5.2356$ \\
		${d_i}\left( 0 \right)$&
		$1.05$ \\
		${\varpi _i}(0)$&
		$10$ \\
		${{\iota _i}}$ &
		$560$ \\
		${\bar d}$&
		$3.0$ \\
		$\rho $, ${\varsigma _i}$ &
		$579.6$ \\
		${{\upsilon _i}}$ &
		$0.00173$ \\
		${\eta _i}$ &
		$0.001$ \\
		${{o_i}}$ &
		$0.002$ \\
		\hline
	\end{tabular}
	\label{tab1}
\end{table}	
\section{SIMULATION RESULTS}
The aim of this section is to exhibit simulation examples that underscore the practicality and effectiveness of the proposed control protocol.
Consider an example of $10$ spacecrafts. According to \cite{modle}, the relative dynamics of the $i$-th spacecraft with respect to the designated reference point can be characterized as
\begin{align}\label{eq65}
	\left\{ \begin{array}{l}
		{{\ddot x}_i} = 2\omega {{\dot y}_i} + \dot \omega {y_i} + {\omega ^2}{x_i} + \frac{{2\mu }}{{{r^3}}}{x_i} + {a_{ix}},\\
		{{\ddot y}_i}  =  - 2\omega {{\dot x}_i} - {{\dot \omega }_i}{x_i} + {\omega ^2}{y_i} - \frac{{2\mu }}{{{r^3}}}{y_i} + a{}_{iy},\\
		{{\ddot z}_i} =  - \frac{{2\mu }}{{{r^3}}}{z_i} + {a_{iz}},
	\end{array} \right.
\end{align}
where $(x_i, y_i, z_i)$ denote the relative position of the $i$-th spacecraft with respect to a reference point. The orbital angular velocity of this reference point is denoted by $\omega$, while $\mu$ represents the Earth's gravitational coefficient. The distance between the reference point and the Earth's barycenter is designated as $r$. Additionally, $a_{ix}$, $a_{iy}$, and $a_{iz}$ correspond to the triaxial control inputs.

\begin{figure}[!t]
	\centerline{\includegraphics[width=0.48\textwidth,height=0.18\textwidth]{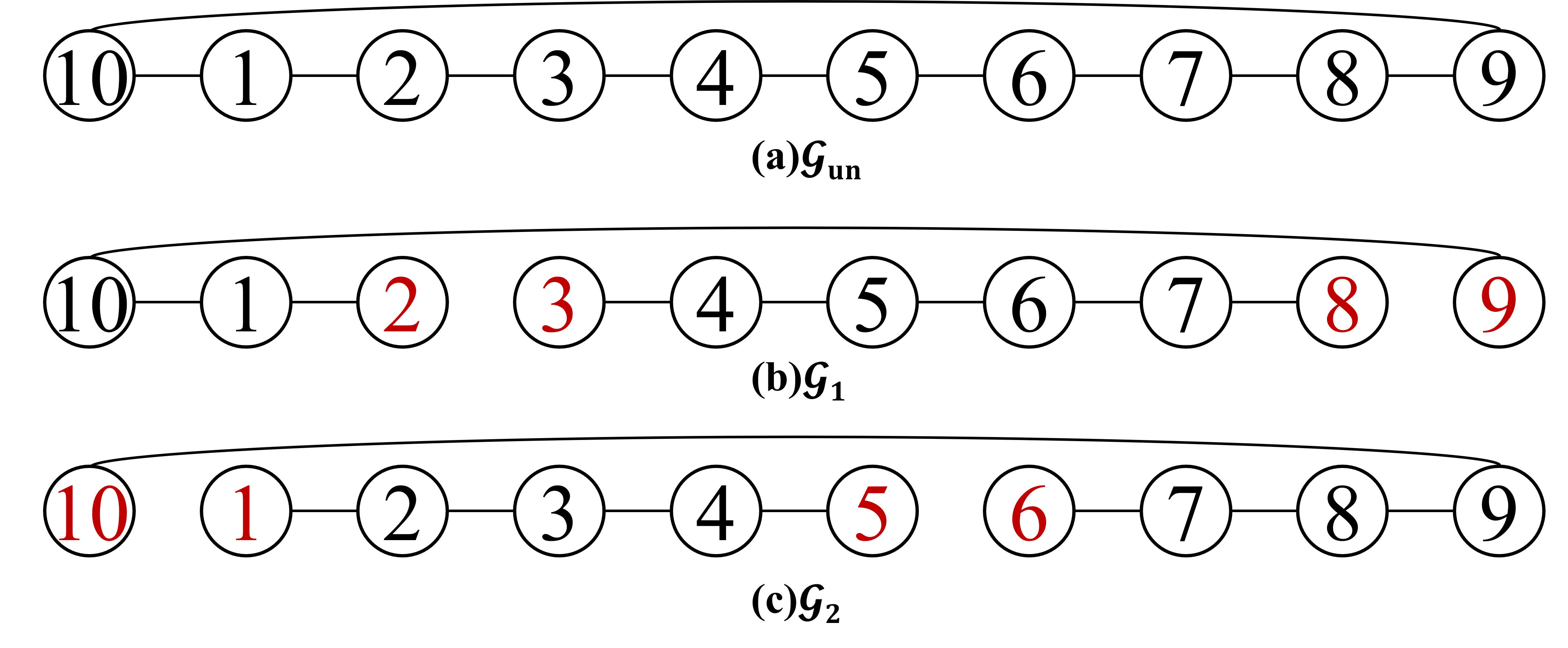}}
	\caption{Communication topologies: (a) union graph of ${{\cal G}_1}$ and ${{\cal G}_2}$; (b) ${{\cal G}_1}$; (c) ${{\cal G}_2}$. In ${{\cal G}_1}$ and ${{\cal G}_2}$, spacecrafts $1$, $2$, $3$, $5$, $6$, $8$, $9$ and $10$ are subject to random communication failure.}
	\label{figga}
\end{figure}

\begin{figure}[!t]
	\centerline{\includegraphics[width=8cm]{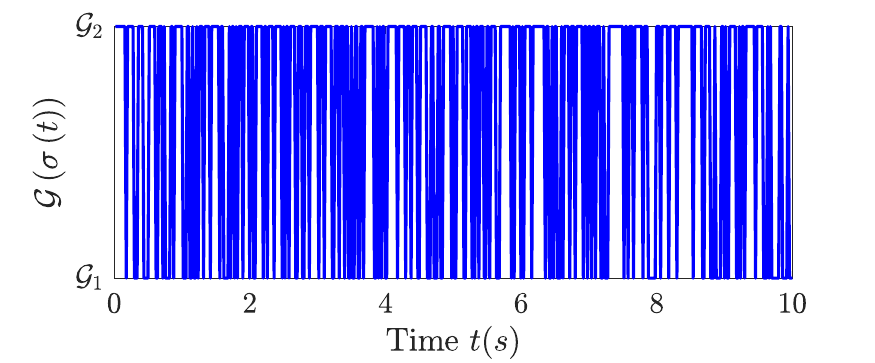}}
	\caption{Illustration depicting the topology switching procedure.}
	\label{figmak3}
\end{figure}

The dynamical model (\ref{eq65}) can be rewritten as state-space model (\ref{eq2}) with
\begin{align}
	&{x_i}\left( t \right) = {\left[ {\begin{array}{*{20}{c}}
				{{x_i}}&{{y_i}}&{{z_i}}&{\begin{array}{*{20}{c}}
						{{{\dot x}_i}}&{{{\dot y}_i}}&{{{\dot z}_i}}
				\end{array}}
		\end{array}} \right]^T} ,\nonumber\\
	&{u_i}\left( t \right) = {\left[ {\begin{array}{*{20}{c}}
				{{a_{ix}}}&{{a_{iy}}}&{{a_{iz}}}
		\end{array}} \right]^T} ,\nonumber
\end{align}
\begin{align}
	A = \left( {\begin{array}{*{20}{c}}
			0&0&0&1&0&0\\
			0&0&0&0&1&0\\
			0&0&0&0&0&1\\
			{{\omega ^2} + \frac{{2\mu }}{{{r^3}}}}&{\dot \omega }&0&0&{2\omega }&0\\
			{ - \dot \omega }&{{\omega ^2} - \frac{\mu }{{{r^3}}}}&0&{ - 2\omega }&0&0\\
			0&0&{{a_{33}}}&0&0&0
	\end{array}} \right) ,\nonumber
\end{align}
\begin{align}
	B = \left( {\begin{array}{*{20}{c}}
			0&0&0\\
			0&0&0\\
			0&0&0\\
			1&0&0\\
			0&1&0\\
			0&0&1
	\end{array}} \right),C = \left( {\begin{array}{*{20}{c}}
			1&0&0&0&0&0\\
			0&1&0&0&0&0\\
			0&0&1&0&0&0
	\end{array}} \right), \nonumber
\end{align}
where the corresponding parameter values are given in Table \ref{table}. We assume that a part of the $10$ spacecrafts are subject to random communication failure, resulting in the topology changes depicted in Figs. \ref{figga}(b) and \ref{figga}(c). The union graph of the communication graphs shown in Figs. \ref{figga}(b) and \ref{figga}(c) is drawn in Fig. \ref{figga}(a), which satisfies Assumption \ref{assumption1}. The diagram depicting the topology switching process utilized in the simulation is presented in Fig. \ref{figmak3}. 
%Assuming that two adjacent spacecrafts can communicate with each other
		\begin{figure}[!t]
	\centerline{\includegraphics[width=8cm]{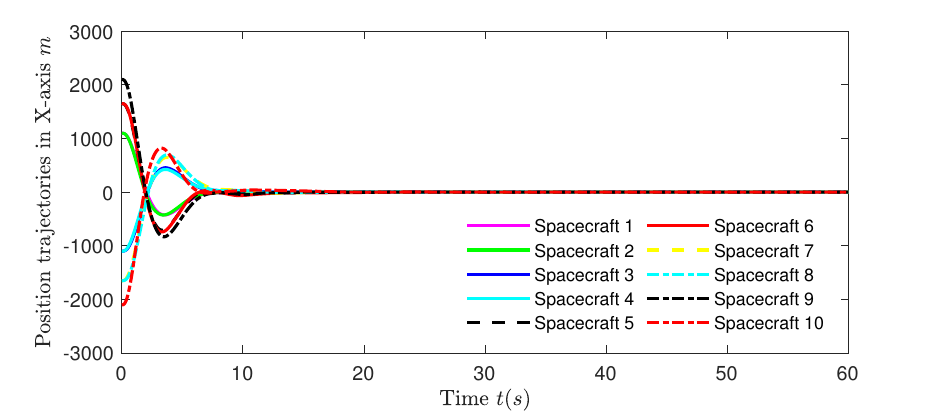}}
	\caption{Trajectories of position in X-axis.}
	\label{fig2}
\end{figure}

\begin{figure}[!t]
	\centerline{\includegraphics[width=8cm]{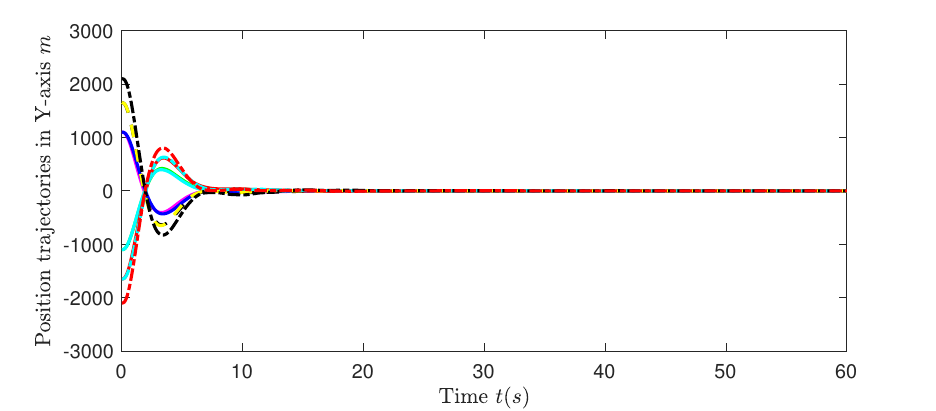}}
	\caption{Trajectories of position in Y-axis.}
	\label{fig3}
\end{figure}

\begin{figure}[!t]
	\centerline{\includegraphics[width=8cm]{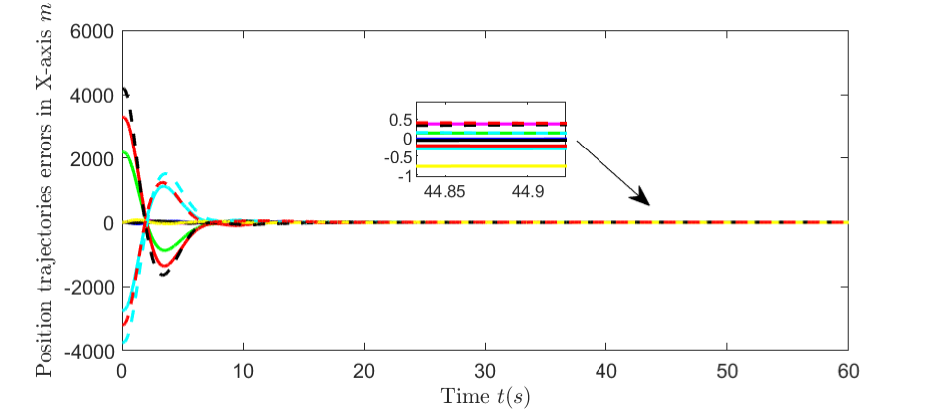}}
	\caption{Trajectories of position errors in X-axis.}
	\label{fig4}
\end{figure}

According to Theorem \ref{thrm3}, we can obtain the gain matrices of the controller and the observer

			\begin{align}
				\scriptsize{ \setlength{\arraycolsep}{1.2pt}
					K = {B^T}P = \left( {\begin{array}{*{20}{c}}
							{0.6016}&{ - 0.0001}&0&{1.7904}&0&0\\
							{0.0001}&{0.6016}&0&0&{1.7904}&0\\
							0&0&{0.6016}&0&0&{1.7904}
					\end{array}} \right) } ,\nonumber
			\end{align}
	{	and
			\begin{align}
				\scriptsize{ \setlength{\arraycolsep}{1.2pt}
					G = {X^{ - 1}}Q = \left( {\begin{array}{*{20}{c}}
							{1.1164}&0&0\\
							0&{1.1164}&0\\
							0&0&{1.1164}\\
							{1.0819}&{0.0001}&0\\
							{0.0001}&{1.0819}&0\\
							0&0&{1.0819}
					\end{array}} \right) }.\nonumber
			\end{align}
			
			Figs. \ref{fig2} and \ref{fig3} show the position trajectories of $10$ spacecrafts along the X-axis and Y-axis, respectively. Figs. \ref{fig4} and \ref{fig5}, on the other hand, present the corresponding position tracking errors, indicating that the steady-state values of position tracking errors are less than $0.5$ $m$. This supports the achievement of the secure consensus for MASs under the distributed observer-based controller. Additionally, Fig. \ref{fig6} depicts the trajectories of the adaptive coupling strengths ${d_i}$, indicating that all the adaptive coupling strengths can converge to a finite value in the end. Fig. \ref{fig7}, finally, displays the event-triggered instants for the $1$-st, $4$-th, $7$-th, and $10$-th spacecrafts. It is noteworthy that there is no Zeno behavior observed in the designed event-triggered mechanism.

			\begin{figure}[!t]
				\centerline{\includegraphics[width=8cm]{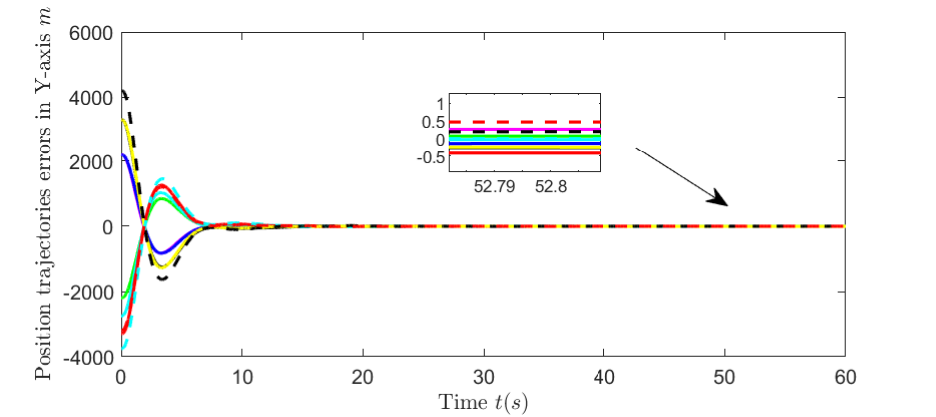}}
				\caption{Trajectories of position errors in Y-axis.}
				\label{fig5}
			\end{figure}
			
			\begin{figure}[!t]
				\centerline{\includegraphics[width=8cm]{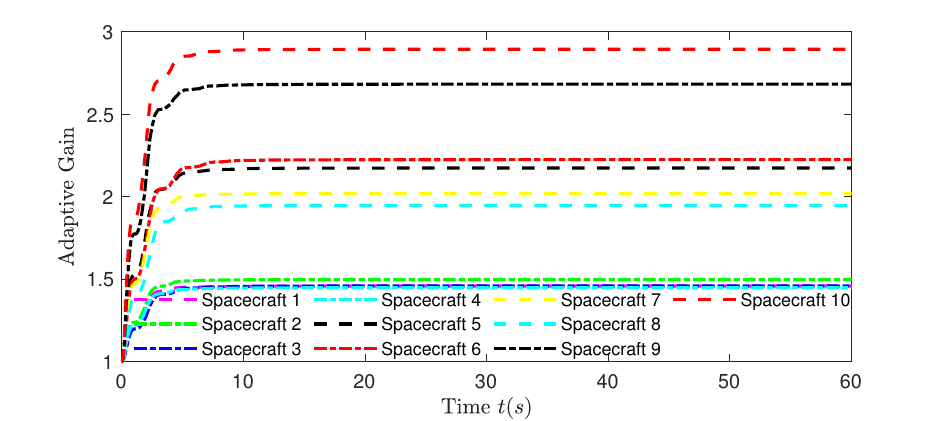}}
				\caption{The time-varying coupling strengths ${d_i}$ of $10$ spacecrafts.}
				\label{fig6}
			\end{figure}
			\begin{figure}[!t]
				\centerline{\includegraphics[width=8cm]{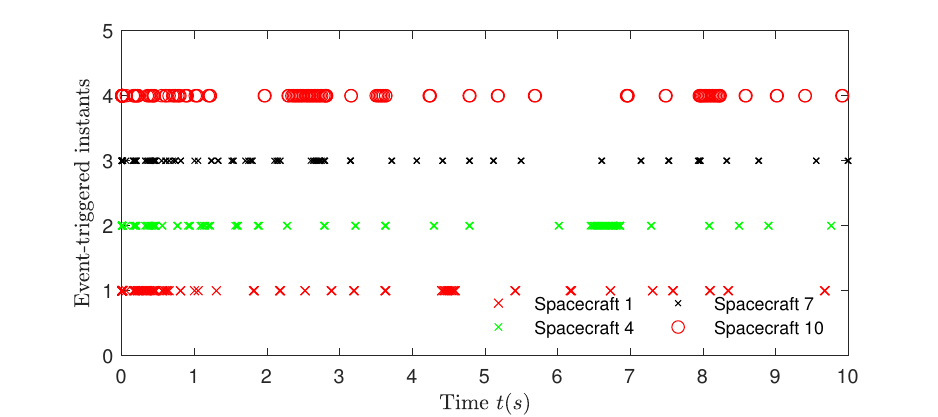}}
				\caption{Event-triggered instants of the $1$-st, $4$-th, $7$-th, $10$-th spacecrafts.}
				\label{fig7}
			\end{figure}

			\begin{figure}[!t]
				\centerline{\includegraphics[width=8cm]{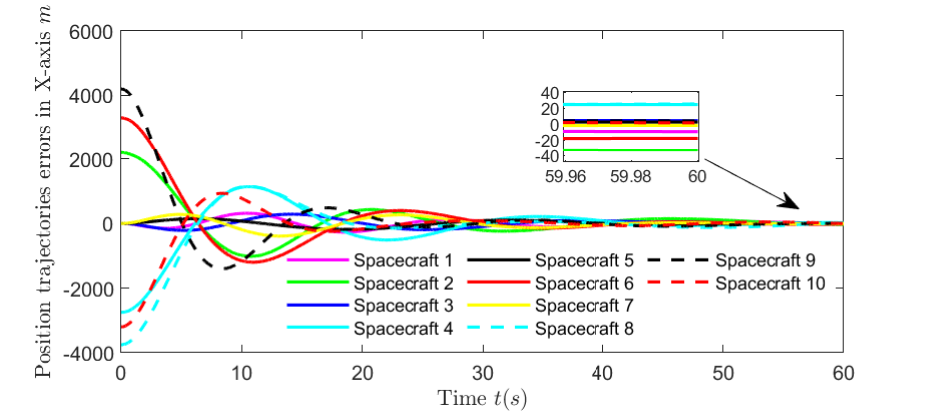}}
				\caption{Trajectories of X-axis position errors under the controller in \cite{simcom}.}
				\label{figcomx1}
			\end{figure}
			
			\begin{figure}[!t]
				\centerline{\includegraphics[width=8cm]{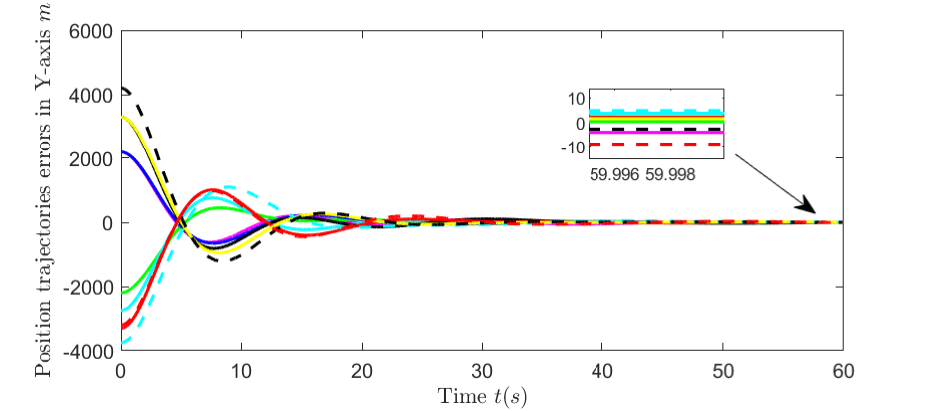}}
				\caption{Trajectories of Y-axis position errors under the controller in \cite{simcom}.}
				\label{figcomx2}
			\end{figure}
			
			\begin{figure}[!t]
				\centerline{\includegraphics[width=8cm]{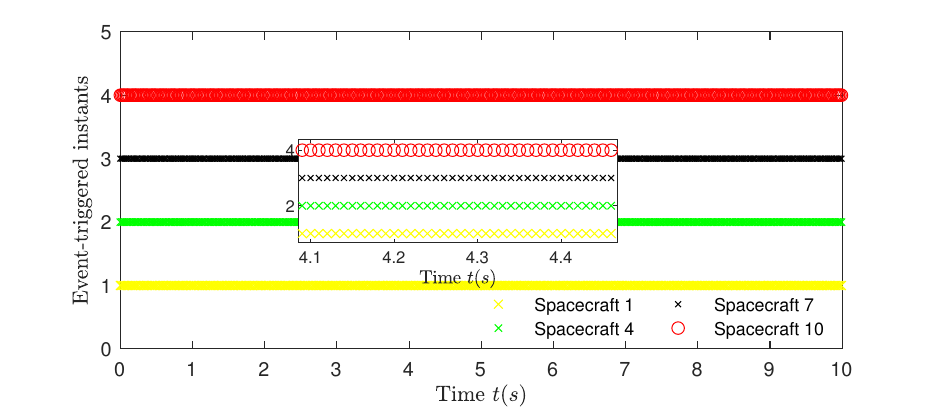}}
				\caption{Event-triggered instants of the $1$-st, $4$-th, $7$-th, $10$-th spacecrafts under the controller in \cite{simcom}.}
				\label{figcomevent}
			\end{figure}
			
			 For comparison, the distributed event-triggered consensus protocol in \cite{simcom} is introduced, which deals with the leaderless consensus control problem for MASs under fixed and directed communication topology. Specially, the consensus protocol and the corresponding event-trigger function are given as
			\begin{IEEEeqnarray}{c}
					{u_i}\left( t \right) =  - \kappa {B^T}P{\xi _i}\left( {t_k^i} \right),\\
			{f_i} = \frac{\kappa }{{{\beta _1}}}e_i^T\left( t \right)PB{B^T}P{e_i}\left( t \right) - {\beta _2}{\gamma _i}\mu \xi _i^T\left( t \right){\xi _i}\left( t \right) - c, \nonumber
			\end{IEEEeqnarray}
			where $\xi_i=\sum_{j=1}^N a_{ij}(x_i-x_j)$, $e_i=\xi_i(t_k^i)-\xi_i(t)$,
			\begin{align}
				\scriptsize{ \setlength{\arraycolsep}{1.2pt}
					{B^T}P = \left( {\begin{array}{*{20}{c}}
							{1.0000}&{ - 0.0001}&0&{1.7321}&0&0\\
							{0.0001}&{1.0000}&0&0&{1.7321}&0\\
							0&0&{1.0000}&0&0&{1.7321}
					\end{array}} \right)} ,\nonumber
			\end{align}
			and choose $\mu  = 1$, ${\beta _1} = 0.8$, ${\beta _2} = 0.35$, $\kappa  = 0.05$ and $c = 10$. Under the same Markovian randomly switching communication topologies and deception attacks, Figs. \ref{figcomx1} and \ref{figcomx2} show the evolution trajectories of position errors for spacecrafts employing the distributed event-triggered controller proposed in \cite{simcom}. These results indicate that all the spacecrafts fail to reach consensus under the controller in \cite{simcom} due to the presence of the deception attacks and the switching communication topologies. Notably, the consensus error value in the steady state is significantly larger than that observed in Figs. \ref{fig4} and \ref{fig5}, potentially rendering it unacceptable in practical applications. Additionally, Fig. \ref{figcomevent} demonstrates the event-triggered instants for some agents utilizing the controller from \cite{simcom}. A comparison of Figs. \ref{fig7} and \ref{figcomevent} reveals that the controller in \cite{simcom} requires a considerably higher amount of communication resources compared to our proposed controller, despite its potentially intolerable consensus performance. These comparisons underscore the effectiveness of our proposed resilient observer-based adaptive event-triggered consensus control approach.
			 
\section{CONCLUSION}
In conclusion, this study addresses the intricate challenges encountered by MASs operating in environments prone to deception attacks and Markovian switching topologies. We proposed a novel observer-based distributed event-triggered secure consensus control protocol based on potentially attacked output information, optimizing network resource utilization. By adopting a distributed approach in designing controller and observer gain matrices, we simplified the system design process. Additionally, we proved that our mechanism precludes Zeno behavior, ensuring system stability. Simulation results validate the superiority of our proposed method compared to existing techniques, highlighting its effectiveness and applicability in event-triggered secure control of MASs.

\section{reference}

\bibliographystyle{IEEEtran}
%\bibliography{IEEEtran}
\bibliography{reference1}

\vfill

\end{document}